\documentclass[12pt]{amsart}

\usepackage[utf8]{inputenc}
\usepackage[T1]{fontenc}
\usepackage[english]{babel}
\usepackage{amsmath,amssymb,amsfonts,amsthm,mathtools}
\mathtoolsset{showonlyrefs}
\usepackage{enumitem}

\usepackage[a4paper, left=2cm, right=2cm, top=3cm, bottom=3cm]{geometry}

\usepackage{xcolor, tikz}
\usepackage[font=small]{caption}
\usepackage{hyperref}

\makeatletter

\IfFileExists{images.tex}{%

  \newcommand{\imagepath}[1]{##1}
}{%
  \newcommand{\imagepath}[1]{images/##1}
}

\newtheorem{theorem}{Theorem}[section]
\newtheorem*{conjecture*}{Conjecture}
\newtheorem{lemma}[theorem]{Lemma}
\newtheorem{proposition}[theorem]{Proposition}

\theoremstyle{definition}

\theoremstyle{remark}
\newtheorem{remark}[theorem]{Remark}

\numberwithin{equation}{section}

\let\epsilon=\varepsilon
\newcommand{\IR}{{\mathbb R}}
\DeclareMathOperator{\BV}{BV}
\newcommand{\avg}[1]{{\overline{#1}}}
\newcommand{\intd}{\,{\operatorname{d}}}
\newcommand{\dd}{{\operatorname{d}}}

\DeclareMathOperator{\dom}{dom}
\newcommand{\defeq}{\mathrel{{\vcentcolon}{=}}}
\newcommand{\eqdef}{\mathrel{{=}{\vcentcolon}}}
\newcommand{\abs}[1]{\mathopen| #1\mathclose|}
\newcommand{\norm}[1]{\mathopen\| #1\mathclose\|}
\newcommand{\intervalcc}[1]{\mathopen[ #1\mathclose]}
\newcommand{\intervalco}[1]{\mathopen[ #1\mathclose[}
\newcommand{\intervaloc}[1]{\mathopen] #1\mathclose]}
\newcommand{\intervaloo}[1]{\mathopen] #1\mathclose[}

\@ifundefined{leqslant}{}{\let\le=\leqslant  \let\leq=\leqslant}
\@ifundefined{geqslant}{}{\let\ge=\geqslant  \let\geq=\geqslant}
\@ifundefined{varnothing}{}{\let\emptyset=\varnothing}

\definecolor{c0}{RGB}{220,76,76}
\definecolor{c1}{RGB}{220,157,76}
\definecolor{c2}{RGB}{76,186,220}
\definecolor{c3}{RGB}{94,172,59}
\definecolor{c4}{RGB}{64,85,255}
\definecolor{c5}{RGB}{116,76,220}

\definecolor{c0}{RGB}{118,42,131}
\definecolor{c1}{RGB}{220,157,76}
\definecolor{c2}{RGB}{27,120,55}
\definecolor{c3}{RGB}{97,120,255}
\definecolor{c4}{RGB}{127,191,123}
\definecolor{c5}{RGB}{175,141,195}
\definecolor{c6}{RGB}{162,184,156}

\definecolor{c0}{RGB}{27,158,119}
\definecolor{c1}{RGB}{217,95,2}
\definecolor{c2}{RGB}{117,112,179}
\definecolor{c3}{RGB}{231,41,138}
\definecolor{c4}{RGB}{102,166,30}
\definecolor{c5}{RGB}{230,171,2}
\definecolor{c6}{RGB}{166,118,29}

\definecolor{lightred}{rgb}{1,0.87,0.87}
\definecolor{lightgreen}{rgb}{0.8,1,0.8}
\hypersetup{%
  bookmarksnumbered, linkcolor=blue, citecolor=red,
  bookmarks, breaklinks,
  linkbordercolor=lightred,
  citebordercolor=lightgreen,
}

\makeatother

\begin{document}

\title[Uniqueness, non-degeneracy, and exact multiplicity of positive solutions]{Uniqueness, non-degeneracy, \\ and exact multiplicity of positive solutions \\ for superlinear elliptic problems}

\author[G.~Feltrin]{Guglielmo Feltrin}

\address{
Department of Mathematics, Computer Science and Physics, University of Udine\\
Via delle Scienze 206, 33100 Udine, Italy}

\email{guglielmo.feltrin@uniud.it}

\author[C.~Troestler]{Christophe Troestler}

\address{
Department of Mathematics, University of Mons\\
Place du Parc 20, B-7000 Mons, Belgium}

\email{christophe.troestler@umons.ac.be}



\date{}

\dedicatory{}

\begin{abstract}
In this paper, we focus our attention on the positive solutions to second-order nonlinear ordinary differential equations of the form $u''+q(t)g(u)=0$, where $q$ is a sign-changing weight and $g$ is a superlinear function.
We exploit the classical shooting approach and the comparison theorem
to present non-degeneracy and exact multiplicity results for positive
solutions. This completes the multiplicity results obtained by Feltrin
and Zanolin. Numerical examples and some related open problems
are also discussed.
\end{abstract}

\maketitle

\section{Introduction}
\label{section-1}

The paper investigates the positive solutions to the Dirichlet boundary value problem
\begin{equation}\label{pb-intro}
\begin{cases}
\, -\Delta u = \omega(x) g(u), &\text{in $\Omega$,}
\\
\, u=0, &\text{on $\partial\Omega$,}
\end{cases}
\end{equation}
where $\Omega\subseteq\mathbb{R}^{N}$ is a bounded domain, $\omega$ is a sign-changing weight, and $g$ is a function with a superlinear growth, i.e., $g(u)\sim u^{p}$ with $p>1$.
This kind of problems is usually called ``superlinear indefinite'', a terminology that has been popularized starting with \cite{HeKa-80}. 
As usual, a solution $u$ of \eqref{pb-intro} is said to be positive if $u>0$ in $\Omega$.

In the last thirty years a great quantity of existence and
multiplicity results for positive solutions to indefinite problems,
both in the ODE and PDE cases, have been obtained using different
techniques, such as topological and variational methods, see for instance
\cite{AlTa-96,AmLG-98,BoGoHa-05,Bu-76,CaDaPa-02,DFGoUb-03,LG-00,ObOm-06},
and also \cite{BoFe-23,Fe-18book} for a quite complete panorama of the
research done in this framework.
An important motivation for the analysis of the superlinear indefinite problem \eqref{pb-intro} comes from the search for stationary solutions to parabolic equations arising in different frameworks, such as reaction-diffusion processes and population dynamics models (see \cite{Ac-09, LG-16} and the references therein).

Our work lies on a line of research initiated by the work of
L\'{o}pez-G\'{o}mez \cite{LG-97,LG-00} concerning the analysis of
the number of positive solutions depending on the nodal behaviour of
$\omega$.
More precisely, L\'{o}pez-G\'{o}mez and his collaborators investigate the existence of positive solutions of problems of the form
\begin{equation*}
\begin{cases}
\, -\Delta u = \lambda u + \omega(x) |u|^{p-1} u, &\text{in $\Omega$,}
\\
\, u=0, &\text{on $\partial\Omega$,}
\end{cases}
\end{equation*}
where $\lambda$ is regarded as a continuation parameter and $p>1$.
In \cite{GRLG-00} G\'{o}mez-Re\~{n}asco and L\'{o}pez-G\'{o}mez conjectured the existence of $2^{m}-1$ positive solutions, whenever $\omega$ has $m$ intervals where it is positive separated by intervals where it is negative and $\lambda$ is sufficiently negative; see also \cite{FeLG-22}.
This conjecture was recently shown to be true~\cite{FeLGS25}.
A related line of research was developed by Gaudenzi, Habets and
Zanolin in \cite{GaHaZa-03}. In this paper, the authors focus on weight functions $\omega$
depending on a real parameter $\mu$ which plays the role of strengthening or
weakening its positive and negative parts. More precisely, $\omega$ is expressed in this manner
\begin{equation*}
  \omega = h^{+} - \mu h^{-},
\end{equation*}
where $\mu>0$ is a parameter, $h \colon \mathopen{[}a,b\mathclose{]}
\to \mathbb{R}$, with $h^{\pm} \defeq \abs{h}/2 \pm h/2$ the
positive and negative parts of $h$.
This special form of $\omega$ was already considered by Bandle, Pozio, Tesei in \cite[Section~2.3]{BaPoTe-87} and by L\'{o}pez-G\'{o}mez in \cite{LG-00}.

In this framework, the analogous of the conjecture in \cite{GRLG-00}
is the following multiplicity result proved in \cite{FeZa-15jde} (see
Theorem~\ref{th-FeZa} for a precise statement), which is the starting point of our research.

\begin{theorem}\label{th-intro-00}
  There exists $\mu^{*}> 0$ such that for every $\mu>\mu^{*}$ the
  Dirichlet problem
  \begin{equation}\label{pb-conj-intro}
    \begin{cases}
      \, u'' + \bigl{(} h^{+}(t) - \mu h^{-}(t) \bigr{)} |u|^{p-1}u = 0,
      \qquad p>1,
      \\
      \, u(a) = 0 = u(b),
    \end{cases}
  \end{equation}
  possesses at least $2^{m}-1$ positive solutions, whenever $h$ has
  $m$ intervals where it is positive separated by intervals where it
  is negative.
\end{theorem}

In \cite{GaHaZa-03}, Gaudenzi, Habets and Zanolin show the existence of three
positive solutions to \eqref{pb-conj-intro} when $h$ has two
positivity intervals separated by a negativity one, that is $m=2$, and
$\mu$ is large. This was later improved by the same authors in
\cite{GaHaZa-04} for $m=3$. Both papers are based on a very precise
phase-plane analysis and a shooting technique (the possibility of
using the shooting method in the general case $m \geq 2$ has been
recently asserted in \cite{BoDaPa-17}, see Remark~2.2 therein).

In \cite{FeZa-15jde}, exploiting a topological degree approach, Feltrin and Zanolin
proved the multiplicity result
for every positive integer $m$ and a general superlinear function $g$. Many other analogous multiplicity results for positive solutions followed in different context: supersublinear nonlinearities, logistic-type nonlinearities, periodic and Neumann boundary conditions, $\phi$-Laplacian operators, etc. (cf.~\cite{BoFeSo-20,BoFeZa-18,FeGi-20,FeSo-18,FeSoTe-22,FeZa-17}).

When analyzing these multiplicity results a natural question arise: Is the number of solutions given by the theorem
(and its variants) optimal?
In other words, one should understand whether  there are examples of problems of the form \eqref{pb-conj-intro} admitting exactly $2^{m}-1$ positive solutions for $\mu$ large.

The main result of the present manuscript, as far as we know for the first time in literature, gives a positive answer to the above question and states the following. For simplicity in the exposition, we now present it in a special case (see the general result in Theorem~\ref{th-exact}).

\begin{theorem}\label{th-intro}
  Let $k\in\mathbb{N}$ with $k\geq 2$. Let
  $h(t)=\sin\bigl(k \pi t /(b-a)\bigr)$.  Then, there exists
  $\mu^{**}>0$ such that for every $\mu>\mu^{**}$ the Dirichlet
  problem \eqref{pb-conj-intro} admits exactly $2^{m}-1$ positive
  solutions, where $m$ is the integer part of $(k+1)/2$.
  These solutions are non-degenerate.
\end{theorem}

Roughly speaking, Theorem~\ref{th-intro} can be considered as a
continuation result, in the sense that we are going to explain in the
following.   We start by remarking that the existence of \textit{at
  least} $2^{m}-1$ positive solutions is given by \cite{FeZa-15jde}
and so we have to prove that there are \textit{at most} $2^{m}-1$
solutions. From \cite{FeZa-15jde}, we observe that the number
$2^{m}-1$ comes from the possibility of prescribing, for a positive
solution, the behavior in each interval $I$ where $h$ is positive
between two possible ones: either the solution is ``small'' on $I$ or
it is ``large'' on $I$ (the solution which is ``small'' in every
interval of positivity is identically zero,
thus one have to subtract $1$ from the total).

The first step towards establishing Theorem~\ref{th-intro}
is to analyze the behaviour of each positive solution
of \eqref{pb-conj-intro} when $\mu\to+\infty$. In particular, we prove
that the limit profile of the solution is zero in the intervals where
$h$ is negative and solves the Dirichlet problem in each interval of
positivity (see Proposition~\ref{lem-muinfty}).
The second step is to prove that the number of limit profiles (as $\mu\to+\infty$) is exactly $2^{m}-1$.
The last step is a ``continuation from infinity'' to guarantee that, for $\mu>0$ sufficiently large, $2^{m}-1$ remains the exact number of positive solutions of \eqref{pb-conj-intro}.

The determination of the number of limit profiles relies on
uniqueness criteria for positive solutions to
\begin{equation}\label{eq-intro}
  \begin{cases}
  \,  u'' + f(t,u) = 0,
  \\
  \,  u(a) = 0 = u(b),
  \end{cases}
\end{equation}
where $f\colon\mathopen{[}a,b\mathclose{]}\times\mathbb{R}\to\mathbb{R}$.
The classical multiplicity result proposed by Moore and Nehari in
\cite{MoNe-59} (they consider $f(t,u)=q(t)u^{3}$ with $q$ piecewise
constant and non-negative) shows that the problem of uniqueness can be
of great complexity even in very simple situations (even in the case
of smooth functions, see Fig.~\ref{fig:non-uniqueness}). Looking at the
literature, one can notice that the question of uniqueness of
positive solutions to \eqref{eq-intro} has received a considerable attention
when the nonlinear term $f$ is non-negative, see for instance
\cite{Co-72,Da-95,ErTa-97,Kw-90,Mo-62,Wo-98} for the ODE case and also
\cite{ADCJT-14,LG-16,NiNu-85} for the PDE case, while the few available
results when $f$ is sign-changing are contained in
\cite{BaPoTe-88,BoFeZa-21,BrTo-10,BrHe-90,HaZa-16ampa,Na-16}.
The approach we adopt in this paper is based on a method introduced by
Kolodner \cite{Ko-55} and subsequently extended by Coffman
\cite{Co-67,Co-72}.  This method is reminiscent of the shooting one
since it consists in reducing the boundary value problem for the
linearization of \eqref{eq-intro} to an initial value problem and
in tracking the point of intersection of the solution with the
$t$-axis (see also \cite{Kw-90} for a survey and the references
therein for a list of the main applications).

To reach our goal however, the uniqueness of the limit profile with a
given number of bumps is not enough.  Its non-degeneracy is also
essential for the continuation argument to work.  In our case, though,
the limit profile does \emph{not} solve an ODE (the limit profile is
not twice differentiable!) and so the meaning of
its non-degeneracy is not clear.  To circumvent this difficulty for
the one-dimensional case of~\eqref{pb-intro}, namely
\begin{equation}
  \label{eq-intro-ag}
  \begin{cases}
  \,  u'' + \bigl(h^+(t) - \mu h^-(t)\bigr) g(u) = 0,
  \\
  \,  u(a) = 0 = u(b),
  \end{cases}
\end{equation}
where $h$ changes sign, we have to work separately on each
sub-interval where $h$ is positive and negative.  To ``compose'' the
non-degeneracy results on each sub-interval (see the proof of
Theorem~\ref{th-exact}), the classical uniqueness results mentioned
above need to be extended to precisely track the sign of the solution to
the linearized equation (see
Propositions~\ref{moroney-improved}
and~\ref{moroney-improved-uniform}).

For the reader's convenience, we end this introduction by listing
the hypotheses used.  In the sequel, we assume that
$h \colon \intervalcc{a,b} \to \IR$ is an $L^{1}$-function
(it will be specified explicitly when we need $h$ to be more regular)
and $g \colon \intervalco{0,+\infty} \to \intervalco{0,+\infty}$ is a
$\mathcal{C}^{1}$-function such that
\begin{equation}
  \label{eq:g*}
  g(0)=0, \qquad g(s)>0,
  \quad \text{for all } s \in \intervaloo{0,+\infty}.
  \tag{$g_{+}$}
\end{equation}
Moreover, we also make use of the following superlinear growth conditions
\begin{equation}
\label{eq:gs}
g'(s) > \dfrac{g(s)}{s}, \quad \text{for all $s\in\mathopen{]}0,+\infty\mathclose{[}$,}
\tag{$g_s$}
\end{equation}
\begin{equation}
  \label{eq:g0}
  g'(0) = \lim_{s\to0^{+}} \dfrac{g(s)}{s} = 0,
  \tag{$g_{0}$}
\end{equation}
and
\begin{equation}
  \label{eq:g-infty}
  \lim_{s\to +\infty} \dfrac{g(s)}{s} = +\infty.
  \tag{$g_{\infty}$}
\end{equation}

\smallskip

The paper is organized as follows.
In Section~\ref{section-3}, exploiting the classical Kolodner--Coffman technique, we refine some uniqueness results and present some technical estimates. Those results are then exploited in Section~\ref{section-4} to give the proof of the main exact multiplicity result, based on the study of limit profiles of positive solutions to \eqref{eq-intro-ag}.
At last, in Section~\ref{sec:numerics}, we present some numerical
experiments to graphically illustrate the main theorem and to shed
some light on the behavior of the branches of solutions as $\mu$ moves
away from~$+\infty$; moreover, some related open questions are discussed.

\section{Uniqueness}
\label{section-3}

In this section, we adapt and extend the classical Kolodner--Coffman
method~\cite{Co-67,Co-72,Ko-55,Kw-90} to have more information on the
linearized equation at a positive solution to the Dirichlet boundary
value problem
\begin{equation}\label{eq:edo-main}
  \begin{cases}
  \,  u'' + q(t) g(u) = 0,  \\
  \,  u(a) = 0 = u(b),
  \end{cases}
\end{equation}
where
$q \in L^1(\intervalcc{a,b}, \intervalco{0, +\infty})$,
$q \not\equiv 0$, and
$g \in \mathcal{C}^1(\intervalco{0, +\infty}, \intervalco{0,+\infty})$
satisfies \eqref{eq:g*} and \eqref{eq:gs}.  Thanks
to \eqref{eq:g*}, we can extend $g$ to a continuous function on the
whole real line by setting $g(u) = 0$ whenever $u \in
\intervaloc{-\infty, 0}$.
For $\alpha \in \intervalco{0, +\infty}$, let $u(\cdot; \alpha)$ be the
unique maximal solution to the differential equation
in~\eqref{eq:edo-main} satisfying the initial conditions
\begin{equation*}
  u(a; \alpha) = 0, \qquad
  u'(a; \alpha) = \alpha.
\end{equation*}
Given that the solution $u(\cdot; \alpha)$ is concave when it is
positive and is affine when it is negative, $u(\cdot; \alpha)$ exists
on $\intervalcc{a,b}$.

When $\alpha > 0$, let us also denote $B(\alpha)$ the first
$t \in \intervaloc{a,b}$ such that $u(t; \alpha) = 0$ (if such a time
$t$ exists). 
Let $\dom B \subseteq \intervaloo{0, +\infty}$ be the
domain of~$B$. 
Notice that if $\alpha \in \dom B$, then $u'(B(\alpha); \alpha) < 0$
and, since $g\equiv 0$ for negative real numbers, $u(s; \alpha) < 0$
and $u'(s; \alpha) = u'(B(\alpha); \alpha)$ for every $s \in
\intervaloc{B(\alpha), b}$.

A first criteria ensuring uniqueness of positive solutions of \eqref{eq:edo-main} is the following.

\begin{lemma}\label{lem-delta}
  Under the above assumptions, if
  \begin{equation}\label{eq-wBa}
    \forall \alpha \in\dom B, \quad
    B(\alpha) = b \ \Rightarrow\
    \partial_{\alpha} u (b;\alpha) < 0,
  \end{equation}
  then problem \eqref{eq:edo-main} has at most one positive solution.
\end{lemma}

\begin{proof}
  We divide the proof into two steps.

  \smallskip
  \noindent
  \textit{Step~1. The map $B$ is continuously differentiable.}

  First of all, it is classical (see e.g.~\cite[Chapter~V]{Ha-82} or
  \cite{RoMa-80}) that $u(\cdot;\alpha)$ is continuously
  differentiable with respect to~$\alpha$.  Moreover, we
  have that
  \begin{equation*}
    u\bigl(B(\alpha); \alpha\bigr) = 0,
    \quad
    u'\bigl(B(\alpha); \alpha\bigr) < 0,
    \quad
    \text{for all } \alpha \in \dom B.
  \end{equation*}
  From an application of the implicit function theorem, we deduce that
  $B$ is continuously differentiable
  (from one side if $\alpha$ lies on the boundary of its domain) and
  \begin{equation}\label{eq:dB}
    \partial_{\alpha} B(\alpha)
    = - \dfrac{
      \partial_{\alpha} u\bigl(B(\alpha); \alpha \bigr)}{
      u'\bigl(B(\alpha); \alpha \bigr)}.
  \end{equation}

  \smallskip
  \noindent
  \textit{Step~2. There exists at most one $\alpha \in \dom B$ such
    that $B(\alpha) = b$.}

  Assume that \eqref{eq:edo-main} has a positive solution i.e., there
  exists $\alpha_0 > 0$ such that $\alpha_0 \in \dom B$ and
  $B(\alpha_0) = b$.  We claim that
  $\intervalco{\alpha_0, +\infty} \subseteq \dom B$
  and that for all $\alpha > \alpha_0$,
  $B(\alpha) \in \intervaloo{0, b}$.  Let
  $\intervaloo{\alpha_0, \alpha_1}$ be the largest connected set such
  that
  \begin{equation}
    \label{eq:Bstays<b}
    \forall \alpha \in \intervaloo{\alpha_0, \alpha_1},\quad
    \alpha \in \dom B \text{ and }
    B(\alpha) < b.
  \end{equation}
  Clearly, \eqref{eq-wBa}--\eqref{eq:dB} implies that
  $\intervaloo{\alpha_0, \alpha_1} \ne \emptyset$. If
  $\alpha_1 = +\infty$, the claim is proved.  If not, let us consider
  $t_1 \in \intervalcc{a, b}$ a limit point of $B(\alpha)$ as
  $\alpha \to \alpha_1$.  The continuity of
  $(t,\alpha) \mapsto u(t; \alpha)$ implies that
  $u(t_1; \alpha_1) = 0$ and that
  $u(\cdot; \alpha_1) \geq 0$ on $\intervalcc{a, t_1}$.
  One cannot have $t_1 = a$ because otherwise
  (tracking the maximums on $\intervalcc{a, B(\alpha)}$) that
  would imply $\alpha_1 = u'(a; \alpha_1) = 0$. Therefore, as
  all roots of $u(\cdot; \alpha_1)$ must be simple,
  $u(\cdot; \alpha_1) > 0$ on $\intervaloo{a, t_1}$ and
  $t_1 = B(\alpha_1)$.  If $t_1 = b$, \eqref{eq:Bstays<b}
  implies that $\partial_\alpha B(\alpha_1) \geq 0$ which
  contradicts~\eqref{eq-wBa}.
  If $t_1 \in \intervaloo{a, b}$, the intermediate value theorem
  implies that $\alpha_1$ is in the interior of $\dom B$ and the
  continuity of $B$ implies that $B < b$ on a neighborhood of
  $\alpha_1$, contradicting the maximality of $\alpha_1$.

  Note that the previous argument
  also forbids any $\alpha_1 < \alpha_0$ such that
  $B(\alpha_1) = b$ as we can replay it with $\alpha_0$ and
  $\alpha_1$ swapped. The claim is thus proved.
\end{proof}

\begin{remark}\label{rem-3.2}
  Notice that \eqref{eq:gs} was not used in the above proof.  Moreover,
  Lemma~\ref{lem-delta} remains valid for other boundary conditions
  \emph{mutatis mutandis}.  For example for the boundary conditions
  $u'(a) = 0 = u(b)$ (or $u(a) = 0 = u'(b)$ as we can swap $a$ and
  $b$), the function $u(\cdot; \alpha)$ is defined as the
  unique solution to
  \begin{equation}
    \label{eq:init-val2}
    \begin{cases}
    \, u'' + q(t) g(u) = 0,\\
    \,  u(a; \alpha) = \alpha, \quad
      u'(a; \alpha) = 0,
    \end{cases}
  \end{equation}
  and $B(\alpha)$ is, as before, the first root of
  $u(\cdot; \alpha)$ in $\intervaloc{a, b}$.  A criterion for
  uniqueness is again~\eqref{eq-wBa}.
\end{remark}

Our aim now is to provide conditions for the applicability of
Lemma~\ref{lem-delta}, namely sufficient conditions that guarantee
$\partial_{\alpha}u(b; \alpha) < 0$, for every
$\alpha \in \intervalco{0,+\infty}$ yielding a solution.
First of all, we recall the following useful version of the well
known Sturm's Comparison Theorem~\cite{Ha-82, Zettl-05}.

\begin{theorem}[Sturm's comparison Theorem]\label{th-Sturm1}
  Let $\omega_{1},\omega_{2} \colon \intervalcc{c,d} \to \mathbb{R}$
  be $L^{1}$-functions with $\omega_{1}(t)\leq \omega_{2}(t)$ in
  $\intervalcc{c,d}$ with a strict inequality on a set of positive
  measure.  Let $u_{1}$ and $u_{2}$ be nontrivial solutions to
  \begin{equation*}
    u_{i}''+\omega_{i}(t)u_{i}=0
    \quad \text{on $\mathopen]c,d\mathclose[$,}
    \quad i=1,2.
  \end{equation*}
  Assume that one of the following conditions holds
  \begin{itemize}
  \item $u_1(c) = 0 = u_1(d)$;
  \item $u'_1(c) = 0 = u_1(d)$ and $u'_2(c) = 0$.
  \end{itemize}
  Then, $u_{2}$ must vanish at some point in $\intervaloo{c, d}$.
\end{theorem}

Then, we have the following.

\begin{lemma}\label{lem-1zero}
  Suppose the assumptions at the beginning of this section hold.
  Then, for all $\alpha > 0$ such
  that $\alpha \in \dom B$, the function
  $\partial_\alpha u(\cdot; \alpha)$ has at least one zero in
  $\intervaloo{a, B(\alpha)}$.
\end{lemma}

\begin{proof}
  Fix $\alpha > 0$ with $\alpha\in \dom B$, and let $\hat b \defeq B(\alpha)$.
  We apply the Sturm's Comparison Theorem~\ref{th-Sturm1} to
  problem~\eqref{eq:edo-main} written in the following form
  \begin{equation*}
    u'' + \dfrac{q(t) g(u)}{u} u = 0,
    \qquad
    u(a) = 0 = u(\hat b),
  \end{equation*}
  and the equation satisfied by
  $w \defeq \partial_\alpha u(\cdot; \alpha)$, namely
  \begin{equation*}
    w'' + q(t) g'(u(t; \alpha)) w = 0.
  \end{equation*}
  Exploiting \eqref{eq:gs}, we deduce that $w$ has at least one zero
  between the two zeros $a$ and $\hat b$ of $u(\cdot; \alpha)$. The
  lemma is proved.
\end{proof}

\begin{remark}\label{rem-3.3}
  As before, Lemma~\ref{lem-1zero} remains valid \textit{mutatis
    mutandis} for the boundary conditions $u'(a) = 0 = u(b)$ (and
  $u(a) = 0 = u'(b)$, swapping $a$ and $b$).
\end{remark}

Given the previous lemma, to establish that
$\partial_{\alpha}u(b; \alpha) < 0$, it is sufficient to show that
$\partial_{\alpha}u(\cdot; \alpha)$ cannot have more than one zero
in $\intervaloc{a, b}$. 
Below we will refine some classical criteria for this. 
These require to linearize
equation~\eqref{eq:edo-main}; thus a stronger regularity on
$q$ is needed, namely that $q$ is of bounded variation.
Recall that the space $\BV(I)$ of functions of bounded variation on an
interval $I$ is defined by~\cite[Definition~7.1]{Leoni17}:
\begin{equation*}
  \BV(I)
  = \{ u \in L^1(I) \mid u' \text{ is a finite signed Radon measure} \}.
\end{equation*}
Here $u'$ denotes the distributional derivative of~$u$.  Each
$u \in \BV(I)$ can be written as the difference of two bounded nondecreasing
functions and thus admits representatives that are discontinuous at an
at most countable number of
points~\cite[Theorem~7.3 and Theorem~2.36]{Leoni17}.
Therefore $\BV(I) \subseteq L^\infty(I)$.
Moreover, the left and right
limits of such representatives exist at any point $x$ in the closure
of $I$ and are
independent of the representative \cite[Theorem~7.3 and
  Theorem~2.17]{Leoni17}.  These limits will respectively be denoted
by $u(x-)$ and $u(x+)$.  A special representative of $u$,
called the \emph{precise representative} and denoted by
$\avg{u}$, averages the left and right limits:
\begin{equation*}
  \avg{u}(x) := \tfrac{1}{2} \bigl(u(x-) + u(x+) \bigr),
  \qquad
  \text{for any } x \in I.
\end{equation*}
The following version of the Fundamental Theorem of Calculus holds
\cite[Theorem~6.25]{Leoni17}: for all $a, b \in I$ with $a \le b$,
\begin{equation}
  \label{eq:ftc}
  \int_{\intervaloc{a,b}} u'(\dd x)
  = u'\bigl( \intervaloc{a,b} \bigr)
  = u(b+) - u(a+).
\end{equation}
If one integrates over the interval $\intervaloo{a,b}$, $u(b+)$ must
be replaced with $u(b-)$.
Finally we also need a generalization of the Leibniz rule.

\begin{proposition}[Leibniz rule]
  \label{prop:Leibniz}
  Let $I \subseteq \IR$ be an interval and $u, v \in \BV(I)$.  Then
  $uv \in BV(I)$ and
  \begin{equation}
    \label{eq:Leibniz-rule}
    (uv)' = u' \, \avg{v} + \avg{u} \, v'.
  \end{equation}
\end{proposition}

A proof of Proposition~\ref{prop:Leibniz} is given
in~\cite[pp.~189--191]{Volpert-Hudjaev85}.  We offer an alternative
elementary proof
for the one-dimensional case for the reader convenience.

\begin{proof}
  The fact that $uv \in BV(I)$ is asserted on page~45
  of~\cite{Leoni17}.
  The right hand side of~\eqref{eq:Leibniz-rule} defines a measure
  because $u'$ and $v'$ are Radon measures and
  $\avg{u}$ and $\avg{v}$ are bounded and continuous except
  possibly at an at most countable number of points.
  To establish~\eqref{eq:Leibniz-rule}, it
  suffices to prove that the measures on both sides
  of~\eqref{eq:Leibniz-rule} coincide on intervals of the form
  $\intervaloc{a,b} \subseteq I$ \cite[Corollary~B.16]{Leoni17}.
  Thanks to~\eqref{eq:ftc},
  $(uv)'\bigl(\intervaloc{a,b}\bigr) = u(b+) v(b+) - u(a+) v(a+)$.
  From \eqref{eq:ftc}, we deduce that
  \begin{equation*}
    \forall x \in \intervaloc{a,b},\qquad
    \avg{v}(x) = v(a+) + \tfrac{1}{2}
    \Bigl( \int_{\intervaloo{a,x}} + \int_{\intervaloc{a,x}} \Bigr)
    v'(\dd y)
  \end{equation*}
  and similarly for $u$.  Therefore
  \begin{equation*}
    \int_{\intervaloc{a,b}} \avg{v}(x) u'(\dd x)
    = v(a+) \bigl(u(b+) - u(a+)\bigr)
    + \int_{\intervaloc{a,b}}  \tfrac{1}{2}
    \Bigl( \int_{\intervaloo{a,x}} + \int_{\intervaloc{a,x}} \Bigr)
    v'(\dd y) u'(\dd x)  .
  \end{equation*}
  A similar computation followed by a permutation of the integrals
  yields
  \begin{equation*}
    \int_{\intervaloc{a,b}} \avg{u}(y) v'(\dd y)
    = u(a+) \bigl(v(b+) - v(a+)\bigr)
    + \int_{\intervaloc{a,b}}  \tfrac{1}{2}
    \Bigl( \int_{\intervaloc{x,b}} + \int_{\intervalcc{x,b}} \Bigr)
    v'(\dd y) u'(\dd x)  .
  \end{equation*}
  Summing the last two integrals gives
  \begin{align*}
    (u' \, \avg{v} + \avg{u} \, v')\bigl(\intervaloc{a,b}\bigr)
    \hspace*{-6em}
    \\
    &= v(a+) \bigl(u(b+) - u(a+)\bigr)
      + u(a+) \bigl(v(b+) - v(a+)\bigr)
      + \int_{\intervaloc{a,b}} \int_{\intervaloc{a,b}} v'(\dd y) u'(\dd x)
    \\
    &= v(a+) \bigl(u(b+) - u(a+)\bigr)
      + u(a+) \bigl(v(b+) - v(a+)\bigr)
      + \bigl(v(b+) - v(a+)\bigr) \bigl(u(b+) - u(a+)\bigr)
    \\
    &= (uv)'\bigl(\intervaloc{a,b}\bigr)
    .
    \qedhere
  \end{align*}
\end{proof}

\medskip

Notice that assuming $q \in \BV(\intervalcc{a,b})$ allows $q$ to be
piecewise constant as well as piecewise affine.  These simple cases
are already not included in the smoothness assumptions of the original
version of the next results.  Furthermore, let us highlight that the
following
proof differs from the one of Moroney~\cite{Mo-62} and  shows that this
result is a consequence of the property~\eqref{eq-wBa}.

\begin{theorem}[Moroney]\label{th-moroney}
  Suppose the assumptions at the beginning of this section hold.
  Assume further that
  $q\in\BV(\intervalcc{a,b})$ is
  non-increasing on $\intervalcc{a,b}$.  Then, \eqref{eq-wBa} holds
  for the map $u$ defined by \eqref{eq:init-val2}.  In particular, the
  differential equation in \eqref{eq:edo-main} with boundary
  conditions $u'(a) = 0 = u(b)$ has at most one positive solution.
\end{theorem}

\begin{proof}
  Let $u(\cdot; \alpha)$ be the solution
  to~\eqref{eq:init-val2}.  Let $\alpha > 0$ be such that
  $u(t; \alpha) > 0$ for all $t \in \intervaloo{a,b}$ and vanishes
  for $t=b$.  Let $w(t) \defeq \partial_\alpha u(t; \alpha)$.  From
  Lemma~\ref{lem-delta}, we know that there is at most one positive
  solution if $w(b) < 0$ for all such $\alpha$'s.  From
  Lemma~\ref{lem-1zero}, we know that $w$ has a (simple) root in
  $\intervaloo{a,b}$.  If we show that $w$ has no other root in
  $\intervaloo{a,b}$, we are done.  This will be done following the
  ideas of Sturm's Separation Theorem~\cite[Theorem~2.6.2]{Zettl-05}.

  Suppose on the contrary that there are two roots of $w$ in
  $\intervaloc{a,b}$.  Thus there exists $c < d$ such that
  \begin{equation}
    \label{eq:w-two-roots}
    w(c) = 0 = w(d), \quad
    w < 0 \text{ on } \intervaloo{c,d}, \quad
    \text{and}\quad
    w'(c) < 0 < w'(d).
  \end{equation}
  Recall that $w$ satisfies $w'' + q(t) g'(u) w = 0$ where
  $u(\cdot) \defeq u(\cdot; \alpha)$.  Now consider
  $u'(\cdot) \defeq u'(\cdot; \alpha)$.
  Differentiating the equation of $u$ and using
  Proposition~\ref{prop:Leibniz}, we obtain that $u'$ satisfies
  \begin{equation*}
    u'''  + \avg{q}(t) g'(u) u' = - q'(t) \avg{g(u)}
  \end{equation*}
  in the sense of measures.
  Note that, the function $g(u)$ being continuous,
  $\avg{g(u)} = g(u)$.
  Moreover $\avg{q} = q$ a.e.\ for the Lebesgue measure (as the number
  of discontinuity points of $q$ is at most countable) and
  consequently also for the measure $g'(u) u' \intd x$.
  Thus $u'$ satisfies
  \begin{equation*}
    u'''  + q(t) g'(u) u' = - q'(t) g(u).
  \end{equation*}
  Multiplying the equation for $w$ (see \eqref{eq:linearized-w}) by
  $u'$ and vice versa and subtracting them gives
  \begin{equation}
    \label{eq:q'}
    \dfrac{\mathrm{d}}{\mathrm{d}t}
    \bigl[ u'(t)w'(t) - u''(t)w(t) \bigr]
    = q'(t) g(u) w(t)
  \end{equation}
  in the sense of measures.
  Here we made use of \eqref{eq:Leibniz-rule} as well as the
  identities $\avg{u'} = u'$, $\avg{w} = w$, $\avg{w'} = w'$, because
  these functions are (absolutely) continuous, and
  $\avg{u''} w' = u'' w'$ because $w' \in L^1$ and so the measure
  $w' \intd x$ does not care about the value of $u''$ at its points of
  discontinuity which are at most countable.
  Then integrating from $c$ to $d$ and using \eqref{eq:ftc} yields:
  \begin{equation*}
    u'(d) w'(d) - u'(c) w'(c)
    = \int_{\intervaloc{c,d}} g(u(t)) w(t) \, q'(\mathrm{d} t).
  \end{equation*}
  Since $q' \le 0$ as a measure on $\intervalcc{a,b}$,
  the right-hand side is $\geq 0$.  On
  the other hand, $u'(t) = \int_a^t -q g(u) \leq 0$ so the
  left-hand side is $\leq 0$.  Therefore both are null and
  $0 = u'(c) = \int_a^c -q g(u)$.  This in turn implies that
  $q \equiv 0$ a.e.~on $\intervalcc{a,c}$ because $g(u) > 0$ on
  $\intervaloo{a,b}$ as a consequence of~\eqref{eq:g*}.
  Therefore $w'' \equiv 0$ on $\intervaloo{a,c}$
  which, given the initial conditions of $w$,
  contradicts the fact that $c$ is a root of $w$.
\end{proof}

The next lemma proposes sufficient conditions on the weight $q$
ensuring the symmetry with respect to the middle point of the interval
$\intervalcc{a,b}$ of every positive solution of \eqref{eq:edo-main}
(cf.~\cite{GiNiNi-79} for an analogous result in the PDE setting
under more restrictive regularity assumptions).

\begin{lemma}\label{lem-symmetry}
  Let
  $g \in \mathcal{C}^1(\intervalco{0, +\infty},
  \intervalco{0,+\infty})$ satisfy \eqref{eq:g*}.  Let
  $q\in L^{1}(\intervalcc{a,b})$ be such that $q\ge 0$,
  $q(t) = q(a+b-t)$ a.e.\ on $\intervalcc{a,b}$, and $q$
  non-decreasing on $\intervalcc{a,c}$, where $c \defeq (a+b)/2$. Then,
  every positive solution $u$ of \eqref{eq:edo-main} is symmetric,
  i.e., $u(t)=u(a+b-t)$, for every $t\in\intervalcc{a,b}$,
  and is non-decreasing on $\intervalcc{a,c}$.
\end{lemma}

\begin{figure}[ht]
\centering
\begin{tikzpicture}[x=0.35cm,y=0.35cm]
\draw[->] (-0.5,0) -- (10.5,0);
\draw[dashed, color = gray] (3.5,0) -- (3.5,7);
\draw (0,-0.2) -- (0,0.2);
\draw (0,-0.8) node {$a$};
\draw (5,-0.2) -- (5,0.2);
\draw (5,-0.8) node {$c$};
\draw (10,-0.2) -- (10,0.2);
\draw (10,-0.8) node {$b$};
\draw (3.5,-0.2) -- (3.5,0.2);
\draw (3.5,-0.8) node {$t^{*}$};
\draw (7,-0.2) -- (7,0.2);
\draw (7,-0.8) node {$\tilde{a}$};
\draw [line width=1.2pt, color = black, rounded corners]
(0,0) .. controls (1,4) and (2,6) .. (3.5,6.5) .. controls (5,6) and (8,4) .. (10,0);
\draw (7,5) node {$u$};
\draw [dashed, line width=1.2pt, color = magenta, rounded corners]
(0,0) .. controls (1,4) and (2,6) .. (3.5,6.5) .. controls (5,6) and (6,4) .. (7,0);
\draw (7,2) node {$\tilde{u}$};
\end{tikzpicture}
\hspace{20pt}
\begin{tikzpicture}[x=0.35cm,y=0.35cm]
\draw[->] (-0.5,0) -- (10.5,0);
\draw[dashed, color = gray] (3.5,0) -- (3.5,5);
\draw[dashed, color = gray] (2,4.7) -- (8,4.7);
\draw[dashed, color = gray] (5,0) -- (5,6);
\draw (0,-0.2) -- (0,0.2);
\draw (0,-0.8) node {$a$};
\draw (5,-0.2) -- (5,0.2);
\draw (5,-0.8) node {$c$};
\draw (10,-0.2) -- (10,0.2);
\draw (10,-0.8) node {$b$};
\draw (3.5,-0.2) -- (3.5,0.2);
\draw (3.5,-0.8) node {$t^{*}$};
\draw (7,-0.2) -- (7,0.2);
\draw (7,-0.8) node {$\tilde{a}$};
\draw [line width=1.2pt, color = black, rounded corners]
(0,1) .. controls (1,4) and (3,4.7) .. (5,5)
.. controls (7,4.7) and (8,4) .. (10,1);
\draw (8.1, 4.2) node {$q$};
\draw [dashed, line width=1.2pt, color = magenta, rounded corners]
(0,1) .. controls (1,4) and (3,4.5) .. (3.5,4.7);
\draw [dashed, line width=1.2pt, color = magenta, rounded corners]
(3.5,4.7) .. controls (4,4.5) and (6,4) .. (7,1);
\draw (7.1, 2.2) node {$\tilde{q}$};
\end{tikzpicture}
\caption{Qualitative representation of $\tilde{u}$ and $\tilde{q}$ defined in \eqref{def-tilde-uq}, where $\tilde{a}=2t^{*}-a$.}
\label{fig-tilde}
\end{figure}
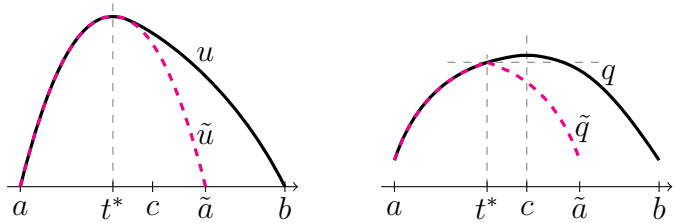

\begin{proof}
Let $u$ be a positive solution of \eqref{eq:edo-main}. 
Since $u$ is concave, the maximum of $u$ in $\intervalcc{a,b}$ is
attained either at a unique maximum point $t^{*}$, or at all the
points of a subinterval $J\subseteq\intervalcc{a,b}$ (if $q\equiv0$ in
$J$), in this latter case let $t^{*} \defeq \min J$. By contradiction, assume that $t^{*}\neq c$.

Due to the symmetry of $q$, $t \mapsto u(a+b-t)$ is another positive solution of \eqref{eq:edo-main}; hence, without loss of generality, we can assume that
$t^{*}\in\intervalco{a,c}$. 
Let $\tilde{u}\colon \intervalcc{a,2t^{*}-a} \to \mathbb{R}$ and $\tilde{q}\colon \intervalcc{a,2t^{*}-a} \to \mathbb{R}$ be defined as follows
\begin{equation}\label{def-tilde-uq}
  \tilde{u}(t) \defeq
\begin{cases}
\, u(t), &\text{if $t\in\intervalcc{a,t^{*}}$,}
\\
\, u(2t^{*}-t), &\text{if $t\in\intervalcc{t^{*},2t^{*}-a}$,}
\end{cases}
\qquad
\tilde{q}(t) \defeq
\begin{cases}
\, q(t), &\text{if $t\in\intervalcc{a,t^{*}}$,}
\\
\, q(2t^{*}-t), &\text{if $t\in\intervalcc{t^{*},2t^{*}-a}$.}
\end{cases}
\end{equation}
See Figure~\ref{fig-tilde} for a graphical representation.
We notice that $\tilde{u}$ and $\tilde{q}$ are symmetric with respect to $t^{*}$ and $\tilde{u}$ satisfies
\begin{equation*}
\begin{cases}
\, \tilde{u}'' + \tilde{q}(t) g(\tilde{u}) = 0,
  & \text{in } \intervalcc{a,2t^{*}-a} \subseteq \intervalcc{a,b},
\\
\, \tilde{u}(a)=0=\tilde{u}(2t^{*}-a).
\end{cases}
\end{equation*}
By the monotonicity and the symmetry of $q$, we observe that
$\tilde{q} \leq q$ in $\intervalcc{a,2t^{*}-a}$.  At last, if
$\tilde{q} < q$ on a set of positive measure, by the Sturm's
Comparison Theorem~\ref{th-Sturm1}, we conclude that $u$ should vanish
in $\intervaloo{a,2t^{*}-a}$, a contradiction.  If instead
$\tilde{q} \equiv q$ a.e.\ on $\intervalcc{a,2t^{*}-a}$, then
$\tilde{u} \equiv u$ on $\intervalcc{a,2t^{*}-a}$
(as they solve the same initial value problem) and so $u(2t^*-a) =
\tilde u(2t^* - a) = 0$, again a contradiction.

The monotonicity of $u$ is a direct consequence of its concavity.
\end{proof}

Combining Lemma~\ref{lem-symmetry} with Theorem~\ref{th-moroney},
we obtain the first criterion for uniqueness to
problem~\eqref{eq:edo-main}.

\begin{theorem}\label{th-moroney-Dirichlet}
  Suppose the assumptions at the beginning of this section hold.
  Let
  $q\in\BV(\intervalcc{a,b})$
  be such that
  $q(t) = q(a+b-t)$ a.e.\ on $\intervalcc{a,b}$ and
  $q$ is non-decreasing on $\intervalcc{a,c}$, where $c = (a + b) / 2$.
  Then, problem
  \eqref{eq:edo-main} has at most one positive solution.
\end{theorem}

Note that this theorem establishes uniqueness without actually proving
that \eqref{eq-wBa} holds.  However, for the next section, this and
actually more information on the solutions of the linearized equation
is necessary.  This is what we do in the next proposition.

\begin{proposition}
  \label{moroney-improved}
  Under the same assumptions as in Theorem~\ref{th-moroney-Dirichlet},
  if $u$ is a positive solution to \eqref{eq:edo-main} and $w$ is a
  nontrivial solution to the linearized equation
  \begin{equation}
    \label{eq:linearized-w}
    \begin{cases}
    \,  w'' + q(t) g'(u) w = 0,\\
    \,  w(a) \ge 0,\ w'(a) \ge 0.
    \end{cases}
  \end{equation}
  Then, $w(b) < 0$ and $w'(b) < 0$.
\end{proposition}

\begin{proof}
  Reasoning as in the proof of Lemma~\ref{lem-1zero}, using Sturm's
  Comparison Theorem~\ref{th-Sturm1}, we show the existence of
  $\xi \in \intervaloo{a, b}$ such that $w(\xi)=0$ and $w(t) > 0$, for
  all $t \in \intervaloo{a, \xi}$.  We claim that
  \begin{equation}\label{eq:claim-xi}
    \xi \in \intervaloo{c, b},
  \end{equation}
  where $c \defeq (a + b) / 2$.  We suppose by contradiction that
  $\xi \in \intervaloc{a, c}$.  Let $\eta \in \intervalco{a, c}$ be
  such that $w'(\eta)=0$.  Such a point $\eta$ exists because of the
  conditions on $w$ at $a$ in~\eqref{eq:linearized-w}.
  Differentiating the equation
  in~\eqref{eq:edo-main}, we obtain, as in the proof of
  Theorem~\ref{th-moroney},
  the following equation for $u'$:
  \begin{equation*}
    u''' + q(t) g'(u) u' = - q'(t) g(u)
  \end{equation*}
  in the sense of measures.
  Then multiplying the equation of $w$ in \eqref{eq:linearized-w} by
  $u'$ and vice versa and subtracting them, we obtain \eqref{eq:q'}.
  An integration from $\eta$ to $\xi$ then yields
  \begin{equation*}
    u'(\xi)w'(\xi) + u''(\eta+)w(\eta)
    = \int_{\intervaloc{\eta,\xi}} g(u(t)) \, w(t) \, q'(\mathrm{d} t).
  \end{equation*}
  The left-hand side is $\le 0$ since $w(\eta) > 0$,
  $u''(\eta+) = -q(\eta+) g(u(\eta)) \le 0$, $w'(\xi) < 0$, and
  $u'(\xi) \ge 0$ (thanks to Lemma~\ref{lem-symmetry}).  The right-hand side is
  $\ge 0$ thanks to the hypothesis $q' \ge 0$ in the sense of
  measures on $\intervalcc{a, c}$.
  Therefore both are null and this implies that $q' \equiv 0$
  on $\intervaloc{\eta, \xi}$ and $q(\eta+) = 0$. 
  Thus $q \equiv 0$ a.e.\ on $\intervalcc{\eta, \xi}$, so
  $w'' \equiv 0$ in $\intervalcc{\eta,\xi}$ and we deduce
  that $w' \equiv w'(\eta) = 0$.  Hence $w$ is constant in
  $\intervalcc{\eta,\xi}$.  This contradiction proves
  the~claim~\eqref{eq:claim-xi}.

  Next, we prove that there are no critical points of $w$ in
  $\intervalcc{\xi, b}$.  By contradiction, let
  $\gamma \in \intervalcc{\xi, b}$ be such that $w'(\gamma)=0$
  (without loss of generality, assume that $\gamma$ is the first
  critical point after $\xi$).  Proceeding as before, we
  integrate~\eqref{eq:q'} from $\xi$ to $\gamma$, obtaining
  \begin{equation*}
    0
    \ge - u''(\gamma+) w(\gamma) - u'(\xi) w'(\xi)
    = \int_{\intervaloc{\xi,\gamma}} g\bigl(u(t)\bigr) \, w(t)
    \, q'(\mathrm{d} t) \ge 0,
  \end{equation*}
  since $w(\gamma)<0$,
  $u''(\gamma+) = -q(\gamma+) g(u(\gamma)) \le 0$, $w'(\xi)<0$,
  $u'(\xi) \le 0$ (due to the symmetry, as before).
  Therefore $q' \equiv 0$ on $\intervaloc{\xi, \gamma}$ and $q(\gamma+) =
  0$.  Thus $q \equiv 0$ a.e.\ on $\intervalcc{\xi, \gamma}$, which in turn
  implies that $w'$ is constant on $\intervalcc{\xi, \gamma}$,
  contradicting the definition of $\eta$.
  In conclusion, no critical
  points of $w$ exist in $\intervalcc{\xi, b}$.  Hence $w' < 0$ in
  $\intervalcc{\xi, b}$ and the proof is complete.
\end{proof}

Let us conclude this section by showing that the previous result holds
uniformly with respect to perturbations of $g'$.

\begin{proposition}
  \label{moroney-improved-uniform}
  Under the assumptions of Theorem~\ref{th-moroney-Dirichlet}, let $u$
  be a positive solution to \eqref{eq:edo-main} and $(G_n)_n$
  be a sequence in $\mathcal{C}(\intervalcc{a,b})$ such that
  $G_n \to g'(u(\cdot))$ uniformly on $\intervalcc{a,b}$.  Further,
  let $(w_n)_n$ be a sequence of solutions to
  \begin{equation}
    \label{eq:linearized-wn}
    \begin{cases}
    \,  w'' + q(t) G_n(t) w = 0,
    \\
    \,  w(a) \geq 0,\ w'(a) \geq 0.
    \end{cases}
  \end{equation}
  Then, there exists $\epsilon > 0$ such that, for
  all $n$ sufficiently large,
  \begin{equation*}
   w_n(b) \le -\epsilon \bigl( w_n(a) + w_n'(a) \bigr)
   \quad\text{and}\quad
   w_n'(b) \le -\epsilon \bigl( w_n(a) + w_n'(a) \bigr). 
  \end{equation*}
\end{proposition}

\begin{proof}
  Let $z_{1,n}$ (resp.\ $z_{2,n}$) be the solution to
  \begin{equation}
    \label{eq:def-z12}
    \begin{cases}
    \,  w'' + q(t) G_n(t) w = 0,\\
    \,  w(a) = 1, \
      w'(a) = 0,
    \end{cases}
    \qquad
    \left(
      \text{resp. }
      \begin{cases}
      \,  w'' + q(t) G_n(t) w = 0,\\
      \,  w(a) = 0, \
        w'(a) = 1
      \end{cases}
    \right) .
  \end{equation}
  We claim that there exists $\epsilon > 0$ such that, for all $n$
  sufficiently large and $j \in \{1,2\}$,
  \begin{equation}
    \label{eq:zjn-sign}
    z_{j,n}(b) \le -\epsilon \text{ and }
    z_{j,n}'(b) \le -\epsilon.
  \end{equation}
  Suppose on the contrary that there exists a $j \in \{1,2\}$ and a
  subsequence, still denoted $(z_{j,n})_n$, such that
  \begin{equation}
    \label{eq:zj-b}
    \liminf_{n \to \infty} z_{j,n} \ge 0
    \quad\text{or}\quad
    \liminf_{n \to \infty} z_{j,n}' \ge 0.
  \end{equation}
  Because $(G_n)_{n}$ is bounded in $\mathcal{C}(\intervalcc{a,b})$ and
  the initial conditions are also bounded, Gr\"onwall's Lemma on
  $z_{j,n}^2 + (z_{j,n}')^2$ implies that the sequence $(z_{j,n})_n$
  is bounded in $\mathcal{C}^1(\intervalcc{a,b})$.
  Now, using again the equation and invoking Ascoli--Arzel\`a Theorem,
  one deduces that  $(z_{j,n})_n$ converges (taking if
  necessary a subsequence) in $\mathcal{C}^{1}$ and one names its
  limit $z_{j,\infty}$.  Thanks to
  $G_n \to g'(u)$, $z_{j,\infty}$ is a
  solution to
  \begin{equation}
    \label{eq:lim-eq-pos}
    \begin{cases}
    \,  z'' + q(t) g'(u(t)) z = 0,\\
    \,  z(a) = 1, \quad
      z'(a) = 0,
    \end{cases}
    \ \text{(if } j = 1 \text{)}
    \qquad\text{or}\quad
      \begin{cases}
      \,  z'' + q(t) g'(u(t)) z = 0,\\
      \,  z(a) = 0, \quad
        z'(a) = 1
      \end{cases}
      \ \text{(if } j = 2 \text{).}
  \end{equation}
  Moreover, \eqref{eq:zj-b} implies that $z_{j,\infty}(b) \ge 0$ or
  $z'_{j,\infty}(b) \ge 0$.
  This contradicts Proposition~\ref{moroney-improved}.

  To conclude the proof, it suffices to notice that
  \begin{equation}
    \label{eq:wn-combili-gen}
    w_n = w_n(a) z_{1,n} + w'_n(a) z_{2,n} 
  \end{equation}
  and therefore \eqref{eq:zjn-sign} implies the claim.
\end{proof}

\section{Exact multiplicity}
\label{section-4}

In this section, we focus our attention on the Dirichlet boundary value problem
\begin{equation}
\label{eq:Dmu}
\begin{cases}
\, u'' + \bigl{(}h^{+}(t)-\mu h^{-}(t) \bigr{)} g(u) = 0, \\
\, u(a) = 0 = u(b),
\end{cases}
\end{equation}
where $\mu$ is a real parameter, $h^{+}$ and $h^{-}$ denote
the positive and, respectively, the negative part of the weight
function $h$, namely $h^{\pm} \defeq |h|/2 \pm h/2$.

In the present section, our starting point is the following result
proposed in \cite{FeZa-15jde} for Dirichlet and mixed boundary
conditions and subsequently in \cite{FeZa-17} for the Neumann and
periodic boundary value problems. In the framework of problem \eqref{eq:Dmu} it reads as follows.

\begin{theorem}[Feltrin and Zanolin, 2015]\label{th-FeZa}
Let $g \colon \mathopen{[}0,+\infty\mathclose{[} \to \mathopen{[}0,+\infty\mathclose{[}$ be a continuous function satisfying \eqref{eq:g*}, \eqref{eq:g0} and~\eqref{eq:g-infty}.
Let $h \colon \mathopen{[}a,b\mathclose{]} \to \mathbb{R}$ be an $L^{1}$-function. Suppose that
\begin{enumerate}[leftmargin=26pt,labelsep=8pt,label=\textup{$(h_{*})$}]
\item there exist $2m+2$ points
\begin{equation*}
a =\tau_{0} \leq \sigma_{1} < \tau_{1} < \cdots < \sigma_{i} < \tau_{i} < \cdots < \sigma_{m} <  \tau_{m}  \leq \sigma_{m+1} = b,
\end{equation*}
such that $h\succ0$ on
$[\sigma_{i},\tau_{i}]$, for every $i\in\{1,\dotsc,m\}$,
and $h\prec0$ on
$[\tau_{i},\sigma_{i+1}]$, for every $i\in\{0,\dotsc,m\}$.
\label{hp-sign-h}
\end{enumerate}
Then, there exists $\mu^{*}>0$ such that for
  every $\mu>\mu^{*}$ problem \eqref{eq:Dmu} has at least
  $2^{m}-1$ positive solutions.
\end{theorem}

The symbol $h\succ 0$ means that $h\geq0$ almost everywhere and
$h\not\equiv0$ on a given interval, and $h \prec 0$ stands for
$-h \succ 0$.  Moreover, for simplicity in the sequel, we set
\begin{equation*}
  I_{i}^{+} \defeq \intervalcc{\sigma_{i},\tau_{i}},
  \quad i\in\{1,\dotsc,m\},
  \qquad\text{and}\qquad
  I_{i}^{-} \defeq \intervalcc{\tau_{i},\sigma_{i+1}},
  \quad i\in\{0,\dotsc,m\},
\end{equation*}
the intervals of positivity and the ones of negativity for the function $h$. Without loss of generality, from now on, we assume that the points $\sigma_{i}$ and $\tau_{i}$ are selected in such a manner that $h\not\equiv0$ on all left neighborhoods of $\sigma_{i}$ and on all right neighborhoods of $\tau_{i}$ (cf.~\cite[Section~5.2]{FeZa-15jde}).

In \cite{FeZa-15jde} the authors introduce a topological approach
based on an extension of the Leray--Schauder degree for locally compact operators on open (possibly
unbounded) sets. 
Along the proof, they first introduce three constants
$0 < r < R$ (with $r$ small and $R$ large) and $\mu^{*}>0$ (sufficiently large) such that for every positive solution $u$ of \eqref{eq:Dmu} with $\mu>\mu^{*}$, it holds that
\begin{equation}\label{prop-r-R}
0 < \max_{t\in I^+_i}|u(t)|<R, 
\quad  \max_{t\in I^+_i}|u(t)| \neq r,  
\quad\text{for every $i \in \{1,\dotsc, m\}$}.
\end{equation}
Next, based on the above mentioned degree theory,
given any nonempty set of indices
$\mathcal{I} \subseteq \{1,\dotsc, m\}$, they prove that there exists at least one
positive solution $u_{\mathcal{I},\mu}$ to \eqref{eq:Dmu} contained in the set
\begin{equation}\label{def-Lambda-I}
  \begin{aligned}
    \Lambda^{\mathcal{I}} \defeq
    \biggl\{\,u\in\mathcal{C}(\mathopen{[}a,b\mathclose{]})\colon
    & r < \max_{t\in I^+_i}|u(t)|<R,\
      i\in\mathcal{I} \text{ and}
    \\
    & \max_{t\in I^+_i}|u(t)| < r, \
      i \in \{1,\dotsc, m\}\setminus\mathcal{I} \; \biggr\},
\end{aligned}
\end{equation}
namely, $u_{\mathcal{I},\mu}$ is ``small'' on $I_{i}^{+}$ for $i\notin\mathcal{I}$ and ``large'' (i.e., $r<\max_{I_{i}^{+}} u_{\mathcal{I},\mu} < R$) if $i\in\mathcal{I}$.
It is worth noting that $u_{\mathcal{I},\mu}$ is concave in each
$I_{i}^{+}$ and convex in each $I_{i}^{-}$, due to the sign condition
on $h$.  As a consequence, $\norm{u}_{\infty} < R$.
Let us also remark that the precise value of $r > 0$ does not matter
as long as it is small enough so that the only non-negative solution
in $\Lambda^\emptyset$ (i.e.,
small on all the intervals $I^+_i$) is the trivial solution
(see \cite[Lemma~2.2]{FeZa-15jde}).

The following result illustrates the convergence of the solutions
for $\mu\to+\infty$.

\begin{proposition}\label{lem-muinfty}
Let $(\mu_{n})_{n} \subseteq \mathopen]\mu^{*}, +\infty\mathclose[$ be such that $\mu_{n} \to +\infty$, and $(u_{\mu_{n}})_{n}$ be a sequence of positive solutions to problem~\eqref{eq:Dmu}.
Then, there exists a continuous function $u_{\infty} \colon \mathopen{[}a,b\mathclose{]} \to \mathopen{[}0,+\infty\mathclose{[}$ such that, going to a subsequence of $(u_{\mu_{n}})_{n}$ if necessary, one has
\begin{equation*}
  \lim_{n \to +\infty} u_{\mu_{n}} = u_{\infty}
\quad \text{uniformly on } \mathopen{[}a,b\mathclose{]}.
\end{equation*}
Moreover, for all $i \in \{0,\dotsc, m\}$, $u_{\infty}\equiv0$ on $I_{i}^{-}$, and the restriction $u_{\infty}|_{I^{+}_{i}} \colon I^{+}_{i}\to\mathopen{[}0,+\infty\mathclose{[}$ is a non-negative solution to
\begin{equation}\label{limit-pb}
\begin{cases}
\, u'' + h^{+}(t) g(u) = 0 , \\
\, u(\sigma_{i}) = 0 = u(\tau_{i}).
\end{cases}
\end{equation}
Furthermore, if in addition
$(u_{\mu_{n}})_{n} \subseteq \Lambda^{\mathcal{I}}$ for some
$\mathcal{I}\subseteq\{1,\dotsc, m\}$ with $\mathcal{I}\neq\emptyset$,
then, for all $i \in \{1,\dotsc, m\} \setminus \mathcal{I}$,
$u_{\infty} \equiv 0$ on $I_{i}^{+}$ and, for all
$i \in \mathcal{I}$, $u_{\infty}|_{I^{+}_{i}}$ is a
positive solution to \eqref{limit-pb}.
\end{proposition}

\begin{proof}
Let $(\mu_{n})_{n} \subseteq \mathopen]\mu^{*}, +\infty\mathclose[$ be such that $\mu_{n} \to +\infty$, and $(u_{\mu_{n}})_{n}$ be a sequence of positive solutions to problem~\eqref{eq:Dmu}.
Then, recalling \eqref{prop-r-R} and \eqref{def-Lambda-I}, we find
\begin{equation*}
  u_{\mu_{n}} \in
  \bigcup_{\mathcal{I}\subseteq\{1,\dotsc, m\},\ \mathcal{I} \ne \emptyset}
  \Lambda^{\mathcal{I}}, \quad \text{for every $n$.}
\end{equation*}
The case $\mathcal{I} = \emptyset$ is excluded because $u_{\mu_n}$ is
nontrivial.
Thus, up to a subsequence, $(u_{\mu_{n}})_{n} \subseteq
\Lambda^{\mathcal{I}}$ for some nonempty $\mathcal{I}\subseteq\{1,\dotsc, m\}$.
From now on, to highlight this property, we denote $(u_{\mathcal{I},\mu_{n}})_{n} \subseteq \Lambda^{\mathcal{I}}$ that subsequence.

The proof borrows some ideas developed in~\cite{BoFeZa-18,FeZa-15jde,FeZa-17}.
Let $\mu>\mu^{*}$ and let $\mathcal{I} \subseteq \{1,\dotsc,m\}$ be
nonempty. Let
$u_{\mathcal{I},\mu}\in\Lambda^{\mathcal{I}}$ be a solution to
problem~\eqref{eq:Dmu}. Then, $\|u_{\mathcal{I},\mu}\|_{\infty} \leq R$ and
\begin{equation}\label{eq:bound-L1-I+}
|u''_{\mathcal{I},\mu}(t)|
\leq h^{+}(t) \max_{s\in\mathopen{[}0,R\mathclose{]}} g(s),
\quad \text{for a.e.~$t\in I^{+}_{i}$,}
\end{equation}
for all $i \in \{1,\dotsc, m\}$.
Next, since $\norm{u_{\mathcal{I},\mu}}_{\infty} \leq R$, the mean
value theorem implies that there exists $t^{*}_{i}\in I^{+}_{i}$ such that
$|u_{\mathcal{I},\mu}'(t^{*}_{i})| = |u_{\mathcal{I},\mu}(\tau_{i})-u_{\mathcal{I},\mu}(\sigma_{i})|/|I^{+}_{i}| \leq R/|I^{+}_{i}|$.
Then, we have
\begin{equation}\label{eq:bound-u'-I+}
  |u'_{\mathcal{I},\mu}(t)|
  = \biggl{|}u'_{\mathcal{I},\mu}(t^{*}_{i}) + \int_{t^{*}_{i}}^{t}
  u''_{\mathcal{I},\mu}(\xi) \,\mathrm{d} \xi \biggr{|}
   \leq \frac{R}{|I^{+}_{i}|}
  + \norm{h}_{L^{1}(I^{+}_{i})} \max_{s\in [0,R]} g(s)
  \eqdef \kappa_{i},
\quad \text{for every $t\in I^{+}_{i}$,}
\end{equation}
for all $i\in\{1,\dotsc,m\}$. As a
consequence of the convexity of $u_{\mathcal{I},\mu}$ on
$I^{-}_{i}$, $|u'_{\mathcal{I},\mu}(\cdot)|$ is bounded on
$\mathopen{[}a,b\mathclose{]}$. Then, via Ascoli--Arzel\`{a} Theorem,
there exists
$u_{\mathcal{I},\infty}\in\mathcal{C}(\mathopen{[}a,b\mathclose{]})$
such that, up to a subsequence,
$u_{\mathcal{I},\mu_{n}} \to u_{\mathcal{I},\infty}$ uniformly
on $\mathopen{[}a,b\mathclose{]}$.

\smallskip
The rest of the proof is divided into three steps.

\smallskip
\noindent
\textit{Step~1.} We are going to prove that $u_{\mathcal{I},\mu}$
tends uniformly to $0$ on all the intervals $I^{-}_{i}$. More
precisely, we claim that for every $\varepsilon$ with $0 < \varepsilon
\leq r$, there exists $\mu^{\star}_{\varepsilon} \geq \mu^{*}$
such that for every $\mu > \mu^{\star}_{\varepsilon}$
and $i\in\{1,\dotsc,m\}$,
we have
$\max_{t\in I^{-}_{i}} u_{\mathcal{I},\mu}(t) < \varepsilon$.

Let $\varepsilon \in \intervaloc{0,r}$.  Let us define
\begin{equation*}
  \delta_{i} \defeq \min \biggl\{\dfrac{|I^{-}_{i}|}{2},
  \dfrac{\varepsilon}{2 \kappa_{i}} \biggr\}
  \quad\text{and}\quad
  \tilde{\delta}_{i} \defeq
  \min \biggl\{\dfrac{|I^{-}_{i}|}{2},
  \dfrac{\varepsilon}{2 \kappa_{i+1}} \biggr\},
\end{equation*}
where $\kappa_i$ and $\kappa_{i+1}$ are defined in
\eqref{eq:bound-u'-I+}, as well as
\begin{align*}
  \mu^{\text{left}}_{i}
  &\defeq \dfrac{R + \kappa_{i} \delta_{i}}{
    \min_{s\in\mathopen{[}\varepsilon/2,R\mathclose{]}} g(s)
    \int_{\tau_{i}}^{\tau_{i} + \delta_{i}} \int_{\tau_{i}}^{t}
    h^{-}(\xi)\,\mathrm{d} \xi \,\mathrm{d} t},
  \\
  \mu^{\text{right}}_{i}
  &\defeq \dfrac{R + \kappa_{i+1} \tilde\delta_{i}}{
    \min_{s\in\mathopen{[}\varepsilon/2,R\mathclose{]}} g(s)
    \int_{\sigma_{i+1} - \tilde\delta_{i}}^{\sigma_{i+1}}
    \int_{t}^{\sigma_{i+1}} h^{-}(\xi)\,\mathrm{d} \xi \,\mathrm{d} t}.
\end{align*}
The denominators are positive because $h^-$ is never identically zero
in right neighborhoods of $\tau_i$ nor in left neighborhoods of
$\sigma_{i+1}$.  We claim that, for
$\mu > \mu^{\star}_{\varepsilon} \defeq
\max_{i=1,\dotsc,m}\{\mu^{\text{left}}_{i}, \mu^{\text{right}}_{i},
\mu^{*}\}$ and $i \in \{1,\dotsc,m\}$, we have
$u_{\mathcal{I},\mu}(\tau_{i}) < \varepsilon$ and
$u_{\mathcal{I},\mu}(\sigma_{i+i}) < \varepsilon$.  If so, by the
convexity of $u_{\mathcal{I},\mu}$ on $I^{-}_{i}$, we immediately have
$u_{\mathcal{I},\mu}(t) < \varepsilon$ for all $t\in I^{-}_{i}$.

Let $\mu > \mu^*_\epsilon$ and $i \in \{1,\dotsc,m\}$.  Suppose on the
contrary that $u_{\mathcal{I},\mu}(\tau_{i}) \ge \varepsilon$ or
$u_{\mathcal{I},\mu}(\sigma_{i+1}) \ge \varepsilon$.

Let us first deal with the case
$u_{\mathcal{I},\mu}(\tau_{i}) \ge \varepsilon$.
By~\eqref{eq:bound-u'-I+},
$u_{\mathcal{I},\mu}'(\tau_{i}) \geq -\kappa_{i}$.  The convexity of
$u_{\mathcal{I},\mu}$ on $I^{-}_{i}$ guarantees that
$u_{\mathcal{I},\mu}'(t) \geq - \kappa_{i}$ for all $t \in I^{-}_{i}$.
From the definition of $\delta_i$, it is clear that
$u_{\mathcal{I},\mu}(t) \ge \varepsilon/2$ for all
$t \in \intervalcc{\tau_{i},\tau_{i} + \delta_{i}}$.  An integration
on
$\intervalcc{\tau_{i},t} \subseteq \intervalcc{\tau_{i},\tau_{i} +
  \delta_{i}}$ yields
\begin{equation*}
  u_{\mathcal{I},\mu}'(t) 
  = u_{\mathcal{I},\mu}'(\tau_{i})
  + \mu \int_{\tau_{i}}^{t} h^{-}(\xi) g(u_{\mathcal{I},\mu}(\xi)) \intd \xi
  \ge - \kappa_{i}
  + \mu \min_{s\in\mathopen{[}\varepsilon/2,R\mathclose{]}} g(s)
  \int_{\tau_{i}}^{t} h^{-}(\xi) \intd \xi.
\end{equation*}
Next, integrating the above inequality on
$\intervalcc{\tau_{i},\tau_{i} + \delta_{i}}$ and using
$\mu > \mu^*_\varepsilon \ge \mu^{\text{left}}_i$, we obtain
\begin{equation*}
  \begin{aligned}
    u_{\mathcal{I},\mu}(\tau_{i} + \delta_{i})
    &= u_{\mathcal{I},\mu}(\tau_{i})
      + \int_{\tau_{i}}^{\tau_{i} + \delta_{i}}
      u_{\mathcal{I},\mu}'(t)\,\mathrm{d} t
    \\
    & \ge - \kappa_{i} \delta_{i}
      + \mu \min_{s\in\mathopen{[}\varepsilon/2,R\mathclose{]}} g(s)
      \int_{\tau_{i}}^{\tau_{i} + \delta_{i}} \int_{\tau_{i}}^{t}
      h^{-}(\xi)\,\mathrm{d} \xi \,\mathrm{d} t
    \\
    &> - \kappa_{i} \delta_{i} + R + \kappa_{i} \delta_{i} = R,
  \end{aligned}
\end{equation*}
a contradiction with the fact that $R$ is a bound
on~$u_{\mathcal{I}, \mu}$.

The case $u_{\mathcal{I},\mu}(\sigma_{i+1}) \ge \varepsilon$ is
similar.  As above, \eqref{eq:bound-u'-I+} implies that
$u_{\mathcal{I},\mu}'(\sigma_{i+1}) \le \kappa_{i+i}$ and so
$u_{\mathcal{I},\mu}'(t) \le \kappa_{i+i}$ for all $t \in I^-_i$.
By definition of $\tilde{\delta}_{i}$, we have
$u_{\mathcal{I},\mu}(t) \ge \varepsilon /2$ for all
$t\in \intervalcc{\sigma_{i+1} - \tilde\delta_{i}, \sigma_{i+1}}$.
Integrating twice and using
$\mu > \mu^*_\varepsilon \ge \mu^{\text{right}}_i$, we obtain the
contradiction
\begin{equation*}
  u_{\mathcal{I},\mu}(\sigma_{i+1} - \tilde\delta_{i})
  \ge - \kappa_{i+1} \tilde\delta_{i}
  + \mu \min_{s\in\intervalcc{\varepsilon/2,R}} g(s)
  \int_{\sigma_{i+1} - \tilde\delta_{i}}^{\sigma_{i+1}}
  \int_{t}^{\sigma_{i+1}}h^{-}(\xi)\,\mathrm{d} \xi \,\mathrm{d} t
  > R.
\end{equation*}

\smallskip
\noindent
\textit{Step~2.} Let $i \in \mathcal{I}$.  We are going to prove that
$(u_{\mathcal{I},\mu_{n}})$ tends uniformly, up to a subsequence, to a
solution of \eqref{limit-pb} in the interval $I^{+}_{i}$.

We already know that $(u_{\mathcal{I},\mu})_{\mu > \mu^*}$ is bounded
on $I^{+}_{i}$ and, due to \eqref{eq:bound-u'-I+},
so is $(u'_{\mathcal{I},\mu})_{\mu > \mu^*}$.  Then thanks
to~\eqref{eq:bound-L1-I+},
$(u''_{\mathcal{I},\mu})_{\mu > \mu^{*}}$ is bounded in $L^1(I^{+}_{i})$ and
equi-integrable. Then, by combining the Dunford--Pettis Theorem with the
Eberlein--\v{S}mulian Theorem, we obtain that (up to a subsequence)
$u''_{\mathcal{I},\mu_{n}}\rightharpoonup v_{\mathcal{I}}$ in
$L^{1}(I^{+}_{i})$, as $n \to +\infty$, for some
$v_{\mathcal{I}}\in L^{1}(I^{+}_{i})$.
Now, by Ascoli--Arzel\`{a} Theorem, we obtain that
$u_{\mathcal{I},\mu_{n}} \to u_{\mathcal{I},\infty}$ in
$\mathcal{C}^{1}(I^{+}_{i})$, as $n \to +\infty$. Therefore,
$u''_{\mathcal{I},\infty} = v_{\mathcal{I}}$ and
$u_{\mathcal{I},\infty} \in W^{1,2}(I^{+}_{i})$.
By Step~1, we already know that $u_{\mathcal{I},\infty}(\sigma_{i}) =
u_{\mathcal{I},\infty}(\tau_{i}) =0$. Finally, passing to the limit
on the weak formulation of
the equation, we have that $u_{\mathcal{I},\infty}|_{I^{+}_{i}}$ is
a non-negative solution to \eqref{limit-pb}.
Moreover,
because $i \in \mathcal{I}$,
$\max_{I^{+}_{i}} u_{\mathcal{I}, \infty} \geq r$, so $u_{\mathcal{I}, \infty}$ is
nontrivial and, exploiting the strong maximum principle, one deduces
that $u_{\mathcal{I},\infty}|_{I^{+}_{i}}$ is a positive solution to
\eqref{limit-pb}.

\smallskip
\noindent
\textit{Step~3.}  Let $i \in \{1, \dotsc, m\} \setminus
\mathcal{I}$. Using hypothesis \eqref{eq:g0} and reasoning as in
Step~2, $u_{\mathcal{I},\mu_n} \to u_{\mathcal{I},\infty}$ and that
$u_{\mathcal{I},\infty}|_{I^+_i}$ is a non-negative solution to
\eqref{limit-pb}
such that $\max_{I^{+}_{i}} u_{\mathcal{I}, \infty} \le r$.
That implies that $u_{\mathcal{I},\infty}\equiv 0$ in the interval
$I^{+}_{i}$ provided that $r >0$ was chosen small enough.
(See~\cite[Proposition~5.4]{BoFeZa-18} for another approach.)
\end{proof}

Proposition~\ref{lem-muinfty} is the key ingredient of our investigation, since it describes the limit profiles (for $\mu\to+\infty$) of the solutions $u_{\mathcal{I},\mu}$.
The goal of the present section is to find conditions such that $2^{m}-1$ is the exact number of solutions to \eqref{eq:Dmu}, knowing that $2^{m}-1$ is the exact number of limit profiles.

We can thus state and prove our main result.

\begin{theorem}\label{th-exact}
Let $g \in \mathcal{C}^1(\mathopen{[}0,+\infty\mathclose{[}, \mathopen{[}0,+\infty\mathclose{[})$ be a function satisfying \eqref{eq:g*}, \eqref{eq:gs}, \eqref{eq:g0} and \eqref{eq:g-infty}. Let $h \colon \mathopen{[}a,b\mathclose{]} \to \mathbb{R}$ be an $L^{1}$-function satisfying \ref{hp-sign-h}. Moreover, we suppose that
\begin{itemize}
\item for every $i\in\{1,\dotsc,m\}$, $h^{+}\in\BV(I_{i}^{+})$,
  $h^{+}$ satisfies
  $h^{+}(t) = h^{+}(\sigma_{i} + \tau_{i} - t)$ for a.e.\
  $t \in \intervalcc{\sigma_{i},\tau_{i}}$ and is
  non-decreasing on $\intervalcc{\sigma_{i},\zeta_{i}}$,
  where $\zeta_{i} \defeq
  (\tau_{i}+\sigma_{i})/2$ is the middle point of $I_{i}^{+}$.
\end{itemize}
Then, there exists $\mu^{**}>0$ such that for every $\mu>\mu^{**}$
problem \eqref{eq:Dmu} has exactly $2^{m}-1$ positive solutions.
These solutions are non-degenerate.
\end{theorem}

\begin{proof}
  By hypotheses \eqref{eq:g*}, \eqref{eq:g0}, \eqref{eq:g-infty}, and
  \ref{hp-sign-h}, Theorem~\ref{th-FeZa} applies and guarantees the
  existence of at least $2^{m}-1$ positive solutions of \eqref{eq:Dmu}
  for every $\mu$ large enough. Each of these solutions belongs to set
  $\Lambda^{\mathcal{I}}$ for every nonempty subset of indices
  $\mathcal{I}\subseteq\{1,\dotsc,m\}$, cf.~\eqref{def-Lambda-I}. We
  are going to prove that, for $\mu$ larger, there exists a unique
  positive non-degenerate solution to problem \eqref{eq:Dmu} in each
  $\Lambda^{\mathcal{I}}$, for every nonempty
  $\mathcal{I}\subseteq\{1,\dotsc,m\}$.

  Let $\mathcal{I} \subseteq \{1,\dotsc, m\}$ be non-empty.  The
  arguments for non-degeneracy and uniqueness are similar.
  If there are degenerate solutions $u_{\mathcal{I}, \mu}$ for $\mu$
  large, that means that there exists a sequence $\mu_n \to +\infty$
  such that the linearized equation
  \begin{equation*}
    \begin{cases}
    \,  w'' + \bigl( h^{+}(t) - \mu_n h^{-}(t) \bigr)
      g'\bigl(u_{\mathcal{I},\mu_n}(t)\bigr) w = 0,
      \\
    \,  w(a) = 0 = w(b)
    \end{cases}
  \end{equation*}
  possesses a nontrivial solution $w_n$.  Analogously, if there are distinct
  solutions in $\Lambda^{\mathcal{I}}$ for $\mu$ large, then there
  exist a sequence $\mu_n \to +\infty$ and
  $u_{\mathcal{I}, \mu_n} \in \Lambda^{\mathcal{I}}$,
  $v_{\mathcal{I}, \mu_n} \in \Lambda^{\mathcal{I}}$ two different
  positive solutions to~\eqref{eq:Dmu}.  Setting
  $w_n \defeq u_{\mathcal{I}, \mu_n} - v_{\mathcal{I}, \mu_n}$ provides a
  nontrivial solution to
  \begin{align}
    \smash{\raisebox{-1.8ex}{$\left\{ \rule{0pt}{3.8ex}\right.$}}
    &w'' + \bigl( h^{+}(t) - \mu_{n} h^{-}(t) \bigr) G_n(t) w = 0,
      \label{eq:almost-linearized}
    \\
    &w(a) = 0 = w(b),
      \label{eq:almost-linearized-bc}
  \end{align}
  where
  \begin{equation}
    \label{eq:def-Gn}
    G_n(t) \defeq
    \int_0^1 g'\bigl(s u_{\mathcal{I}, \mu_n} (t)
    + (1-s) v_{\mathcal{I}, \mu_n}(t) \bigr)   \intd s.
  \end{equation}
  Thus, if non-degeneracy (resp.\ uniqueness) fails for $\mu$ large,
  there exists a sequence $\mu_n \to +\infty$ such that, for all $n$,
  system~\eqref{eq:almost-linearized}--\eqref{eq:almost-linearized-bc}
  with $G_n(t) \defeq g'\bigl(u_{\mathcal{I},\mu_n}(t)\bigr)$ (resp.\ with
  $G_n$ given by~\eqref{eq:def-Gn}) has a nontrivial solution.  We
  will show hereafter that this is not possible.

  Before that, let us collect the properties of $G_n$ (common to
  both cases).  First, Proposition~\ref{lem-muinfty} says that
  $u_{\mathcal{I},\mu_n} \to u_{\mathcal{I},\infty}$ uniformly on
  $\intervalcc{a,b}$.  In case we have two solutions, one also has
  $v_{\mathcal{I},\mu_n} \to v_{\mathcal{I},\infty}$ uniformly on
  $\intervalcc{a,b}$.  Since
  $u_{\mathcal{I},\infty} \equiv 0 \equiv v_{\mathcal{I},\infty}$ on all
  $I_i^-$, $i \in \{0,\dotsc,m\}$, and $I_i^+$ with
  $i \notin \mathcal{I}$ and we are in a situation where we have
  uniqueness of the positive solution of~\eqref{limit-pb}
  on $I_i^+$, $i \in \mathcal{I}$
  (see Theorem~\ref{th-moroney-Dirichlet}), one necessarily has
  $u_{\mathcal{I},\infty} \equiv v_{\mathcal{I},\infty}$ on
  $\intervalcc{a,b}$.  Therefore, $G_n$ satisfies the following
  properties:
  \begin{itemize}
  \item $G_n \in \mathcal{C}(\intervalcc{a,b})$;
  \item $G_n \to g'(u_{\mathcal{I},\infty}(\cdot))$ uniformly on
    $\intervalcc{a,b}$;
  \item for all $t \in \intervalcc{a,b}$, $G_n(t) \ge 0$.
  \end{itemize}
  In the argument below we will also make use of the following function
  \begin{equation*}
    r_n(t) \defeq \frac{w_{n}(t)}{w'_{n}(t)}
    \quad
    \text{defined for } t \in \intervalcc{a,b}
    \text{ such that } w'_{n}(t) \ne 0.
  \end{equation*}
  Using equation~\eqref{eq:almost-linearized}, one deduces that $r_n$
  satisfies the equation
  \begin{equation}
    \label{eq:ratio-edo}
    r_{n}'(t)
    = \dfrac{ (w_{n}'(t))^{2} - w_{n}(t) w_{n}''(t) }{(w_{n}'(t))^{2}}
    = 1 + \bigl(h^{+}(t) - \mu_n h^{-}(t) \bigr)
    G_n(t) \, r_{n}^{2}(t),
    \quad \text{for a.e.~$t\in\dom r_n$.}
  \end{equation}
  Note that, because $w_n$ is a nontrivial solution
  to~\eqref{eq:almost-linearized}, $w_n$ and $w'_n$ cannot vanish at
  the same time and so $\abs{r_n(t)} \to +\infty$ when $t$ approaches
  the boundary of $\dom r_n$.

  Now, let us rule out the existence of nontrivial solutions to
  \eqref{eq:almost-linearized}--\eqref{eq:almost-linearized-bc} by
  examining in turn how boundary signs transfer for each sub-interval.

  \medskip
  \noindent
  \textit{Step 1.} %
  First, let us consider a negativity interval
  $I_{i}^{-} = \intervalcc{\tau_{i},\sigma_{i+1}}$ for some fixed
  $i \in \{0,\dotsc,m\}$.  From equation~\eqref{eq:almost-linearized}
  we deduce that $w_{n}''(t)\le 0$ for a.e.~$t \in I_{i}^{-}$
  such that $w_{n}(t) < 0$ and $w_{n}''(t) \ge 0$ for a.e.~$t \in I_{i}^{-}$
  such that $w_{n}(t) > 0$.  Therefore, $w_{n}$ is
  concave in the interval $I_{i}^{-}$ whenever $w_{n}$ is negative,
  and convex in the interval $I_{i}^{-}$ whenever $w_{n}$ is positive.
  Therefore, we have
  \begin{alignat}{2}
    &\text{if } w_{n}(\tau_{i}) \ge 0 \text{ and } w_{n}'(\tau_{i}) \ge 0, 
    &&\text{then }
       w_{n}(\sigma_{i+1}) \ge w_{n}(\tau_{i}) \text{ and }
       w_{n}'(\sigma_{i+1}) \ge w_{n}'(\tau_{i});
       \label{eq:1.2-2} \\
    &\text{if } w_{n}(\tau_{i}) \le 0 \text{ and } w_{n}'(\tau_{i}) \le 0,
    &\text{ }&\text{then }
               w_{n}(\sigma_{i+1}) \le w_{n}(\tau_{i})    \text{ and }
               w_{n}'(\sigma_{i+1}) \le w_{n}'(\tau_{i}) .
     \label{eq:1.2-1}
  \end{alignat}
  In cases \eqref{eq:1.2-2}--\eqref{eq:1.2-1}, we claim that, for any
  $c > 0$, one also has
  \begin{equation}
    \label{eq:bound-I-}
    \abs{w_n(\tau_i)} \le c \abs{w'_n(\tau_i)}
    \ \Rightarrow\
    \abs{w_n(\sigma_{i+1})} \le \hat c \abs{w'_n(\sigma_{i+1})} ,
  \end{equation}
  where $\hat c \defeq c + \abs{I_i^-}$.  Indeed the premise of
  \eqref{eq:bound-I-} implies that $w'_n(\tau_i) \ne 0$ (as $w_n$ is a
  nontrivial solution).  Thus, noticing that
  \eqref{eq:almost-linearized} implies that $\abs{w'_n}$ is
  non-decreasing on $I^-_i$, $I^-_i \subseteq \dom r_n$. Then
  equation~\eqref{eq:ratio-edo} implies that $r'_n \le 1$ on $I_i^-$
  and, integrating, we obtain
  \begin{equation*}
    r_n(\sigma_{i+1}) \le r_n(\tau_i) + \abs{I_i^-},
  \end{equation*}
  which proves the claim.

  \medskip
  \noindent
  \textit{Step 2.} %
  Let us now deal with intervals
  $I^+_i = \intervalcc{\sigma_i, \tau_i}$ with $i \notin \mathcal{I}$.
  We first claim that, for $n$ large enough and for any $c > 0$,
  \begin{equation}
    \label{eq:bound-I+-notactivated}
    \abs{w_n(\sigma_i)} \le c \abs{w_n'(\sigma_i)}
    \ \Rightarrow\
    \abs{w_n(\tau_i)} \le \hat c \abs{w_n'(\tau_i)}
  \end{equation}
  where $\hat c \defeq c + 2\abs{I^+_i}$.  To prove that, note that the
  premise of the implication reads $\abs{r_n(\sigma_i)} \le c$
  ($w_n'(\sigma_i) \ne 0$ because the premise
  would otherwise imply
  that $w_n$ is the trivial solution).  We will establish that, for
  $n$ sufficiently large,
  \begin{equation}
    \label{eq:bound-I+-notactivated-rn}
    \forall t \in I^+_i,\qquad
    t \in \dom r_n
    \quad\text{and}\quad
    \abs{r_n(t)} \le \hat c = c + 2\abs{I^+_i},
  \end{equation}
  from which \eqref{eq:bound-I+-notactivated} follows.
  Let $\vartheta_{i}$ be such that
  \begin{equation*}
    0 < \vartheta_{i}
    < \dfrac{|I_{i}^{+}|}{\norm{h}_{L^{1}(I_{i}^{+})} \, \hat c^{2}}.
  \end{equation*}
  Given that $G_n \to g'(u_{\mathcal{I},\infty}) = 0$ uniformly on
  $I_i^-$ (remembering~\eqref{eq:g0}), we can consider $n$ be
  sufficiently large so that $G_n(t) \le \vartheta_{i}$ for all
  $t \in I_{i}^{+}$.  The properties in
  \eqref{eq:bound-I+-notactivated-rn} are plainly satisfied whenever
  $t$ is close to $\sigma_i$.  Let
  $\intervalco{\sigma_{i},t^*} \subseteq
  \intervalcc{\sigma_{i},\tau_{i}}$ be the maximal interval where the
  properties in \eqref{eq:bound-I+-notactivated-rn} are valid.  If
  $t^* = \tau_i$ we are done.  Otherwise,
  $\abs{r_{n}(t^*)} = \hat c$ and $t^* \in \dom r_n$
  (recall that $\abs{r_n}$ blows up when $t$ approaches
  the boundary of $\dom r_n$).
  Then \eqref{eq:ratio-edo} yields the following contradiction:
  \begin{align*}
    2 \abs{I_{i}^{+}}
    = \hat c - c
    &\le\abs{r_{n}(t^*)} - \abs{r_{n}(\sigma_{i})}
      \le \abs
      {r_{n}(t^*) - r_{n}(\sigma_{i})}
      = \biggl| \int_{\sigma_{i}}^{t^*} r_{n}'(t) \intd t \biggr|\\
    &= \biggl| \int_{\sigma_{i}}^{t^*}
      1 + h^{+}(t) G_n(t) r_{n}^{2}(t) \intd t \biggr|
      \le \abs{I_{i}^{+}} + \norm{h}_{L^{1}(I_{i}^{+})} \,
      \vartheta_{i} \, \hat c^{2}
      < 2 \abs{I_{i}^{+}} .
  \end{align*}
  
  Now we claim that, for $n$ large (how large depends on the $c$
  that will be chosen later to apply \eqref{eq:bound-I+-notactivated}),
  \begin{alignat}{2}
    &\text{if } w_{n}(\sigma_{i}) \ge 0 \text{ and } w_{n}'(\sigma_{i}) \ge 0, 
    &\text{ }
    &\text{then }
      w_{n}(\tau_{i}) \ge \frac{1}{2} w_n(\sigma_{i}) \text{ and }
      w_{n}'(\tau_{i}) \ge \frac{1}{2} w_{n}'(\sigma_{i});
      \label{eq:1.3-2} \\
    &\text{if } w_{n}(\sigma_{i}) \le 0 \text{ and } w_{n}'(\sigma_{i}) \le 0, 
    &\text{ }&\text{then }
               w_{n}(\tau_{i}) \le \frac{1}{2} w_{n}(\sigma_{i}) \text{ and }
               w_{n}'(\tau_{i}) \le \frac{1}{2} w_{n}'(\sigma_{i}).
               \label{eq:1.3-1}
  \end{alignat}
  Let $z_{1,n}$ (resp.\ $z_{2,n}$) be as in \eqref{eq:def-z12} with
  $\intervalcc{a,b}$ replaced by
  $I^+_i = \intervalcc{\sigma_i, \tau_i}$.  As in that proof, we have
  \eqref{eq:wn-combili-gen} and $(z_{1,n})$ (resp.\ $(z_{2,n})$)
  converges, taking if necessary a subsequence, in
  $\mathcal{C}^{1}(I^+_i)$ to a function $z_{1,\infty}$ (resp.\
  $z_{2,\infty}$) which is a solution to \eqref{eq:lim-eq-pos} with
  $u \equiv u_{\mathcal{I}, \infty} \equiv 0$ on $I^+_i$. Therefore
  $z_{1,\infty}(t) = 1$ and $z_{2,\infty}(t) = t - \sigma_i$.  Thus,
  for $n$ sufficiently large, $z_{1,n}(\tau_i) \ge 1/2$,
  $z_{1,n}'(\tau_i) \ge -1/(4 c)$,
  $z_{2,n}(\tau_i) \ge \abs{I^+_i}/2$, and $z_{2,n}'(\tau_i) \ge 3/4$.

  If $w_n(\sigma_i) \ge 0$ and $w_n'(\sigma_i) \ge 0$, then the
  previous inequalities and \eqref{eq:bound-I+-notactivated} yield
  \begin{align*}
    w_n(\tau_i)
    &= w_n(\sigma_i) z_{1,n}(\tau_i) + w'_n(\sigma_i) z_{2,n}(\tau_i)
      \ge \frac{1}{2} w_n(\sigma_i) + \frac{\abs{I^+_i}}{2} w'_n(\sigma_i)
      \ge \frac{1}{2} w_n(\sigma_i) ,\\[1\jot]
    w_n'(\tau_i)
    &= w_n(\sigma_i) z_{1,n}'(\tau_i) + w'_n(\sigma_i) z_{2,n}'(\tau_i)
      \ge -\frac{1}{4 c} w_n(\sigma_i) + \frac{3}{4} w'_n(\sigma_i)
      \ge \frac{1}{2} w'_n(\sigma_i) .
  \end{align*}
  Assertion~\eqref{eq:1.3-1} is established in a similar way.

  \medskip
  \noindent
  \textit{Step 3.} %
  At last, let us focus on the intervals
  $I^{+}_{i} = \intervalcc{\sigma_i, \tau_i}$ for some
  $i \in \mathcal{I}$.  We claim that, for any $c > 0$, there exists
  $\varepsilon > 0$ (independent of $n$) such that, for $n$ large
  enough,
  \begin{alignat}{2}
    &\text{if } w_{n}(\sigma_{i}) \ge 0,\ w_{n}'(\sigma_{i}) \ge 0
      \text{ and } \abs{w_n(\sigma_i)} \le c \abs{w_n'(\sigma_i)}\\
    &\text{then }
      w_{n}(\tau_{i}) \le -\varepsilon \bigl( w_{n}(\sigma_{i})
      + w_{n}'(\sigma_{i}) \bigr) < 0, \
      w_{n}'(\tau_{i}) \le -\varepsilon \bigl(
      w_{n}(\sigma_{i}) + w_{n}'(\sigma_{i}) \bigr) < 0 ,\\
    &\text{and }
      \abs{w_n(\tau_i)} \le \hat c \abs{w_n'(\tau_i)}
      \text{ for some } \hat c \text{ independent on } n;
      \label{eq:1.1-1}  \\[1\jot]
    &\text{if } w_{n}(\sigma_{i}) \le 0,\ w_{n}'(\sigma_{i}) \le 0
      \text{ and } \abs{w_n(\sigma_i)} \le c \abs{w_n'(\sigma_i)}\\
    &\text{then }
      w_{n}(\tau_{i}) \ge -\varepsilon \bigl(
      w_{n}(\sigma_{i}) + w_{n}'(\sigma_{i}) \bigr) >0, \
      w_{n}'(\tau_{i}) \ge -\varepsilon \bigl(
      w_{n}(\sigma_{i}) + w_{n}'(\sigma_{i}) \bigr) >0 , \\
    &\text{and }
      \abs{w_n(\tau_i)} \le \hat c \abs{w_n'(\tau_i)}
      \text{ for some } \hat c \text{ independent on } n.
      \label{eq:1.1-2}
  \end{alignat}
  Let us prove \eqref{eq:1.1-1}, \eqref{eq:1.1-2} being similar.
  Let $c > 0$ be fixed.
  Proposition~\ref{moroney-improved-uniform} implies that
  $w_n(\tau_i) \le -\epsilon \bigl( w_n(\sigma_i) + w_n'(\sigma_i) \bigr)$ and
  $w_n'(\tau_i) \le -\epsilon \bigl( w_n(\sigma_i) + w_n'(\sigma_i) \bigr)$
  for $n$ large enough and for some $\epsilon > 0$ independent of $n$.
  Given that $w_n$ can be written as a linear combination of $z_{1,n}$
  and $z_{2,n}$ (see \eqref{eq:wn-combili-gen} where
  $\intervalcc{a,b}$ is here $\intervalcc{\sigma_i, \tau_i}$), one has
  \begin{equation*}
    \frac{w_n(\tau_i)}{w_n'(\tau_i)}
    = r_n(\tau_i)
    = \frac{r_n(\sigma_i) z_{1,n}(\tau_i) + z_{2,n}(\tau_i)}{
      r_n(\sigma_i) z_{1,n}'(\tau_i) + z_{2,n}'(\tau_i)}
  \end{equation*}
  with $r_n(\sigma_i) \in \intervalcc{0, c}$ (recall that, as above,
  $w'_n(\sigma_i) \ne 0$).  Note that
  \begin{equation*}
    \frac{\gamma z_{1,n}(\tau_i) + z_{2,n}(\tau_i)}{
      \gamma z_{1,n}'(\tau_i) + z_{2,n}'(\tau_i)}
    \xrightarrow[n \to +\infty]{}
    \frac{\gamma z_{1,\infty}(\tau_i) + z_{2,\infty}(\tau_i)}{
      \gamma  z_{1,\infty}'(\tau_i) + z_{2,\infty}'(\tau_i)}
    \quad\text{uniformly w.r.t. } \gamma \in \intervalcc{0,c},
  \end{equation*}
  where $z_{1,\infty}$ and $z_{2,\infty}$ are defined in the proof of
  Proposition~\ref{moroney-improved-uniform}.  As
  $z_{1,\infty}'(\tau_i) < 0$ and $z_{2,\infty}'(\tau_i) < 0$ by
  Proposition~\ref{moroney-improved}, the limit is bounded
  independently of $\gamma \in \intervalcc{0,c}$ and thus,
  for $n$ large enough, so is $r_n(\tau_i)$.  In other words, there
  exists a constant $\hat c$ independent of $n$ such that, for $n$ large
  enough, $\abs{w_n(\tau_i)} \le \hat c \abs{w_n'(\tau_i)}$.
  Incidentally, note that this shows that the denominator does not
  vanish for $n$ large and so that $r_n(\tau_i)$ is well defined.
  
  \medskip
  \noindent
  \textit{Step 4.} %
  To conclude the proof of the non-degeneracy and uniqueness, let us
  derive a contradiction.  For this, we will consecutively examine
  the sub-intervals $I^-_i$ and $I^+_i$ starting from $a$ and show
  that, for $n$ large, $w_n(b) \ne 0$.

  Without loss of generality, we can suppose that $w_n'(a) = 1$.  If
  there is a first nontrivial negativity interval
  $I^-_0 = \intervalcc{a, \sigma_1}$, \eqref{eq:1.2-2} implies that
  $w_n(\sigma_1) \ge 0$, $w_n'(\sigma_1) \ge 1$, and
  $\abs{w_n(\sigma_1)} \le c_1 \abs{w_n'(\sigma_1)}$ with
  $c_1 \defeq 1 + \abs{I^-_0}$.  If there is no such first interval, then
  $\sigma_1 = \tau_0 = a$ and we have $w_n(\sigma_1) = 0$,
  $w_n'(\sigma_1) = 1$, and
  $\abs{w_n(\sigma_1)} \le c_1 \abs{w_n'(\sigma_1)}$ with $c_1 \defeq 1$.
  Thus, in both cases we have, for $n$ large, that
  \begin{equation*}
    w_n(\sigma_1) \ge 0, \quad
    w_n'(\sigma_1) \ge 1, \quad\text{and}\quad
    \abs{w_n(\sigma_1)} \le c_1 \abs{w_n'(\sigma_1)} .
  \end{equation*}
  Next we have the interval $I^+_1 = \intervalcc{\sigma_1, \tau_1}$.
  If $1 \notin \mathcal{I}$, then \eqref{eq:bound-I+-notactivated} and
  \eqref{eq:1.3-2} 
  imply that, for $n$ possibly larger,
  \begin{equation}
    \label{eq:I1+-nonactive}
    w_n(\tau_1) \ge 0, \quad
    w_n'(\tau_1) \ge \frac{1}{2}, \quad\text{and}\quad
    \abs{w_n(\tau_1)} \le c_2 \abs{w_n'(\tau_1)} ,
  \end{equation}
  with $c_2 \defeq c_1 + 2 \abs{I^+_1}$.
  If otherwise $1 \in \mathcal{I}$, then \eqref{eq:1.1-1} yield
  \begin{equation}
    \label{eq:I1+-active}
    w_n(\tau_1) \le -\epsilon_1, \quad
    w_n'(\tau_1) \le -\epsilon_1, \quad\text{and}\quad
    \abs{w_n(\tau_1)} \le c_2 \abs{w_n'(\tau_1)} ,
  \end{equation}
  for some $\epsilon_1 > 0$ and $c_2 > 0$ (independent of $n$).

  The next interval is $I_1^- = \intervalcc{\tau_1, \sigma_2}$.  If
  \eqref{eq:I1+-nonactive} occurred, then \eqref{eq:1.2-2} and
  \eqref{eq:bound-I-} yield, possibly taking $n$ larger,
  \begin{equation*}
    w_n(\sigma_2) \ge 0, \quad
    w_n'(\sigma_2) \ge \frac{1}{2}, \quad\text{and}\quad
    \abs{w_n(\sigma_2)} \le c_3 \abs{w_n'(\sigma_2)} ,    
  \end{equation*}
  where $c_3 \defeq c_2 + \abs{I_1^-}$.  If \eqref{eq:I1+-active} occurred,
  then \eqref{eq:1.2-1} and \eqref{eq:bound-I-} yield
  \begin{equation*}
    w_n(\sigma_2) \le -\epsilon_1, \quad
    w_n'(\sigma_2) \le -\epsilon_1, \quad\text{and}\quad
    \abs{w_n(\sigma_2)} \le c_3 \abs{w_n'(\sigma_2)} ,    
  \end{equation*}
  where again $c_3 \defeq c_2 + \abs{I_1^-}$.

  Continuing this procedure through the finitely many sub-intervals
  $I^+_i$ and $I^-_i$, we have to go through (at least) one interval
  $I^+_i$ with $i\in\mathcal{I}$ and thus 
  $|w_n(\sigma_{i})| \geq \epsilon >0$, for $n$ large enough and some
  $\epsilon > 0$ (independent of $n$).
  As a consequence, one ends up with the fact that $w_n(b) \ne 0$
  for $n$ large enough. This contradiction concludes the proof.
\end{proof}

\section{Numerical experiments and open problems}
\label{sec:numerics}

\newif\ifverif  \veriftrue
\veriffalse

\newcommand{\graphSol}[5][]{
  \IfFileExists{\imagepath{#3.dat}}{%
    \draw[->] (0,0) -- (1, 0) -- +(2ex, 0) node[below]{$t$};
    \draw[->] (0,0) -- (0, #5) -- +(0, 1.1ex) node[left]{$u$};
    \foreach \x in {0, 1} {
      \draw (\x, 3pt) -- (\x, -3pt) node[below]{$\scriptstyle \x$};
    }
    \foreach \y in {#4} {
      \draw (1pt, \y) -- (-1pt, \y);
      \node at (-4pt, \y) {$\scriptscriptstyle \y$};
    }
    #1
    \begin{scope}
      \draw[#2, line width=1pt] plot file{\imagepath{#3.dat}};
      \ifverif
        \begin{scope}
          \clip (0, -5) rectangle (1, 10);
          \draw[#2, dashed] plot file{/tmp/#3-lin.dat};
        \end{scope}
      \fi
    \end{scope}
  }{}%
}
\newcommand{\graphMu}[4][]{
  \graphSol[{#1}]{c0}{#2-0}{#3}{#4}
  \begin{scope}[yshift=18mm]
    \graphSol[{#1}]{c1}{#2-1}{#3}{#4}
  \end{scope}
  \begin{scope}[yshift=36mm]
    \graphSol[{#1}]{c2}{#2-2}{#3}{#4}
  \end{scope}
}

This section is devoted to the numerical exploration of questions
related to Theorem~\ref{th-exact} as well as discussions about
possible future investigations.
Throughout this section, we will specialize the function $g$ to the
representative case $g(u) = u^3$ and take $\intervalcc{a,b} =
\intervalcc{0,1}$.

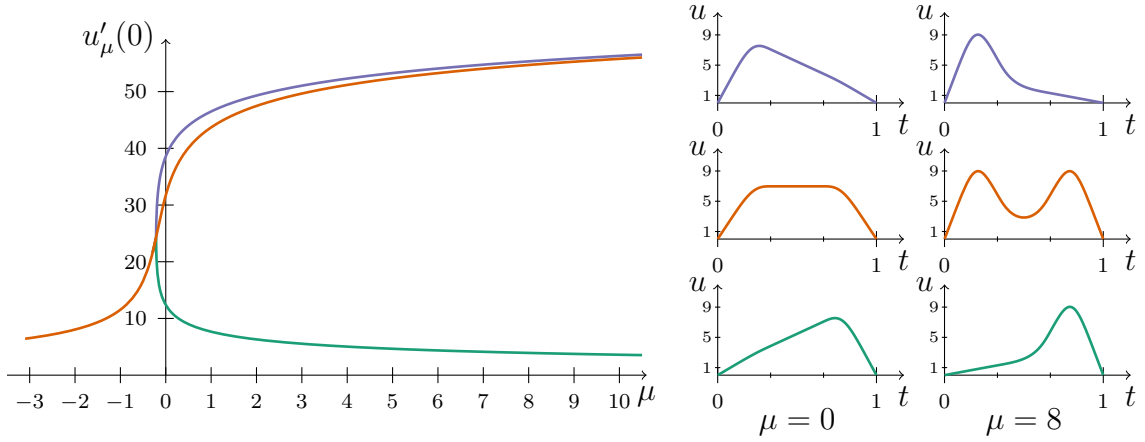
\begin{figure}[htb]
  \centering
  \begin{tikzpicture}[x=6mm, y=0.75mm]
    \draw[->] (-3.5, 0) -- (10.6, 0) node[below]{$\mu$};
    \draw[->] (0, 0) -- (0, 52) -- +(0, 3ex) node[left]{$u'_\mu(0)$};
    \foreach \x in {-3,..., 10} {
      \draw (\x, 3pt) -- (\x, -3pt) node[below]{$\scriptstyle \x$};
    }
    \foreach \y in {10, 20,..., 50} {
      \draw (3pt, \y) -- (-3pt, \y) node[left]{$\scriptstyle \y$};
    }
    \foreach \i in {0, 2, 1}{
      \ifverif
        \draw[color=c\i, line width=0.7pt,
          mark=*, only marks, mark size=0.5pt]
        plot file{\imagepath{branches-3pi-\i.dat}};
        \draw[color=c\i, line width=3pt, opacity=0.3]
        plot file{/tmp/branches-3pi-\i.dat};
      \else
        \draw[color=c\i, line width=1pt]
        plot file{\imagepath{branches-3pi-\i.dat}};
      \fi
    }
    \ifverif
      \draw[line width=0.7pt, color=c5, mark=*, only marks, mark size=0.5pt]
      plot file{\imagepath{branches-3pi-3.dat}};
    \else
      \draw[color=c1, line width=1pt]
      plot file{\imagepath{branches-3pi-3.dat}};
    \fi
    \foreach \x/\m in {73mm/0, 103mm/8}{
      \begin{scope}[xshift=\x, x=21mm, y=1mm]
        \node[below] at (0.5, -3mm) {$\mu = \m$};
        \graphMu[\foreach \x in {0.3333333, 0.6666666}{
            \draw (\x, 1pt) -- (\x, -1pt);
          }]{sol-3-\m}{1, 5, 9}{10}
      \end{scope}
    }
  \end{tikzpicture}
  \caption{The branches of positive solutions to \eqref{eq:Dmu}
    with $h(t) = \sin(3\pi t)$ in $\intervalcc{0,1}$
    and $g(u)=u^{3}$ (on the left) and
    graphs of solutions for some $\mu$ (on the right).}
  \label{fig:bifurc-3pi}
\end{figure}

\begin{figure}[htb]
  \centering
  \newcommand{\MarkSubIntervals}{
    \foreach \x in {0.3333333, 0.6666666}{
      \draw (\x, 1pt) -- (\x, -1pt);
    }
  }
  \begin{tikzpicture}[x=21mm, y=2mm]
    \foreach \m/\x in {-5/0mm, -3/35mm, -1/70mm}{
      \begin{scope}[xshift=\x]
        \node[below] at (0.5, -3mm) {$\mu = \m$};
        \graphSol[\MarkSubIntervals]{c1}{sol-3-\m-0}{1, 4}{4.5}
      \end{scope}
    }
  \end{tikzpicture}
  \caption{Graphs of the unique positive solution to \eqref{eq:Dmu}
    with $h(t) = \sin(3\pi t)$ in $\intervalcc{0,1}$
    and $g(u)=u^{3}$ for various $\mu < 0$.}
  \label{fig:bifurc-3pi-sols}
\end{figure}
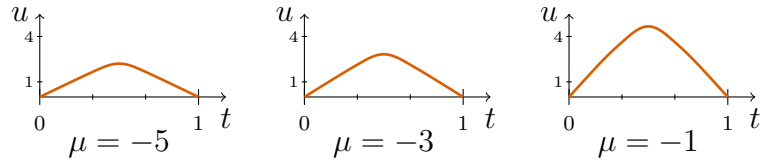

In the previous section, the exact number of positive solutions
of~\eqref{eq:Dmu} for $\mu$ large was determined.  A natural question
is for what range of $\mu$ does this number of solutions persist.  On
Fig.~\ref{fig:bifurc-3pi} and~\ref{fig:bifurc-3pi-sols}, one can see
the branches of positive solutions to~\eqref{eq:Dmu} for
$h(t) = \sin(3\pi t)$ (which corresponds to $m = 2$ in
Theorem~\ref{th-intro}) as well as the graphs of these solutions for
various values of $\mu$.  In this situation, we observe that there
are $3$ positive solutions for all $\mu \ge 0$ and the 3 branches
collapse to a single point when $\mu \approx -0.21$.  Below this
value, a unique positive solution (whence symmetric) exists.  This
bifurcation point can be seen as a symmetry breaking of the ground
state: numerical evidence suggests that,
on the right of the bifurcation point, the value
$\mathcal{A}(u)$ of the action functional
\begin{equation*}
  \mathcal{A}(u)
  \defeq \frac{1}{2} \int_a^b (u')^2 \intd t
  - \int_a^b \bigl(h^+(t) - \mu h^-(t)\bigr) G(u) \intd t,
\end{equation*}
where $G$ is a primitive of $g$, is lower on the ``external'' branches
(the solutions on these branches are symmetric to each other) and
higher on the ``central'' branch made of symmetric solutions.


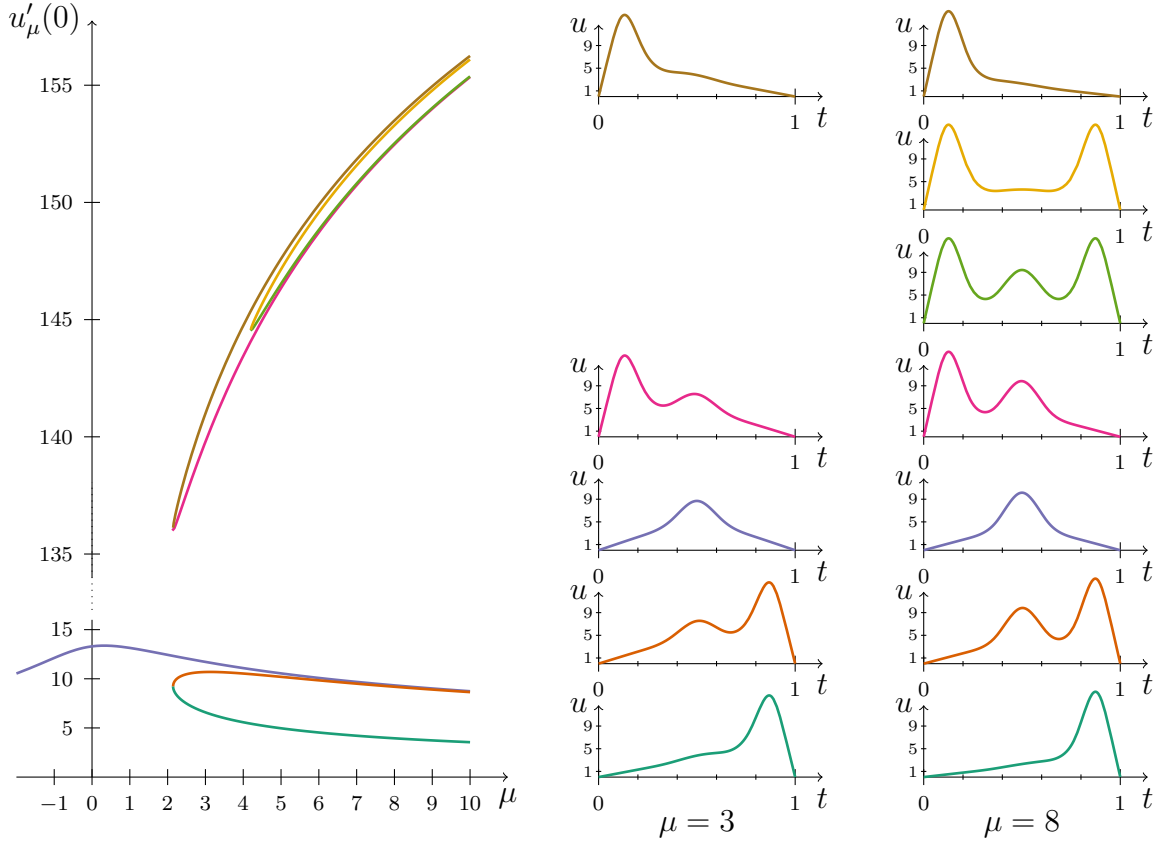
\begin{figure}[htb]
  \centering
  \begin{tikzpicture}[x=5mm, y=1.3mm]
    \draw[->] (-2, 0) -- (11, 0) node[below]{$\mu$};
    \draw (0, 0) -- (0, 16);
    \draw[dotted] (0, 17) -- (0,30);
    \foreach \x in {-1,..., 10} {
      \draw (\x, 3pt) -- (\x, -3pt) node[below]{$\scriptstyle \x$};
    }
    \foreach \y in {5, 10, 15}{
      \draw (3pt, \y) -- (-3pt, \y) node[left]{$\scriptstyle \y$};
    }
    \foreach \i in {0, 2, 1}{
      \ifverif
        \draw[color=c\i, line width=3pt, opacity=0.3]
        plot file{/tmp/branches-5pi-\i.dat};
        \draw[c\i, line width=0.6pt,  mark=*, only marks, mark size=0.5pt]
        plot file{\imagepath{branches-5pi-\i.dat}};
      \else
        \draw[c\i, line width=1pt]
        plot file{\imagepath{branches-5pi-\i.dat}};
      \fi
    }
    \begin{scope}[yshift=-389mm, y=3.1mm]
      \draw[->] (0, 134) -- (0, 156) -- +(0, 3ex) node[left]{$u'_\mu(0)$};
      \foreach \y in {135, 140, 145, 150, 155}{
        \draw (3pt, \y) -- (-3pt, \y) node[left]{$\scriptstyle \y$};
      }
      \foreach \i in {3, 4, 5, 6}{
        \ifverif
          \draw[color=c\i, line width=2.5pt, opacity=0.2]
          plot file{/tmp/branches-5pi-\i.dat};
          \draw[color=c\i, line width=0.6pt,
            mark=*, only marks, mark size=0.5pt]
          plot file{\imagepath{branches-5pi-\i.dat}};
        \else
          \draw[color=c\i, line width=1pt]
          plot file{\imagepath{branches-5pi-\i.dat}};
        \fi
      }
    \end{scope}
    \newcommand{\MarkSubIntervals}{
      \foreach \x in {0.2, 0.4, 0.6, 0.8}{
        \draw (\x, 1pt) -- (\x, -1pt);
      }
    }
    \begin{scope}[xshift=67mm]
      \foreach \i/\c/\y in {0/0/0mm, 1/1/15mm, 2/2/30mm, 3/3/45mm, 4/6/90mm}{
        \begin{scope}[yshift=\y, x=26mm, y=0.75mm]
          \graphSol[\MarkSubIntervals]{c\c}{sol-5-3-\i}{1, 5, 9}{10}
        \end{scope}
      }
      \node[below] at (13mm, -3mm) {$\mu = 3$};
    \end{scope}
    \begin{scope}[xshift=110mm]
      \foreach \i/\y in {0/0mm, 1/15mm, 2/30mm, 3/45mm, 4/60mm,
        5/75mm, 6/90mm}{
        \begin{scope}[yshift=\y, x=26mm, y=0.75mm]
          \graphSol[\MarkSubIntervals]{c\i}{sol-5-8-\i}{1, 5, 9}{10}
        \end{scope}
      }
      \node[below] at (13mm, -3mm) {$\mu = 8$};
    \end{scope}
  \end{tikzpicture}
  \caption{The branches of positive solutions to \eqref{eq:Dmu}
    with $h(t) = \sin(5\pi t)$ in $\intervalcc{0,1}$
    and $g(u)=u^{3}$ (on the left) and
    graphs of solutions for some $\mu$ (on the right).}
  \label{fig:bifurc-5pi}
\end{figure}

\begin{figure}[htb]
  \newcommand{\MarkSubIntervals}{
    \foreach \x in {0.2, 0.4, 0.6, 0.8}{
      \draw (\x, 1pt) -- (\x, -1pt);
    }
  }
  \centering
  \begin{tikzpicture}[x=26mm, y=1.2mm]
    \foreach \x/\m in {0mm/-1, 40mm/0, 80mm/1}{
      \begin{scope}[xshift=\x]
        \node[below] at (0.5, -3mm) {$\mu = \m$};
        \graphSol[\MarkSubIntervals]{c2}{sol-5-\m-0}{1, 5}{6}
      \end{scope}
     }
  \end{tikzpicture}
  \caption{Graphs of the unique positive solution to \eqref{eq:Dmu}
    with $h(t) = \sin(5\pi t)$ in $\intervalcc{0,1}$
    and $g(u)=u^{3}$ for various $\mu \le 1$.}
  \label{fig:bifurc-5pi-unique}
\end{figure}

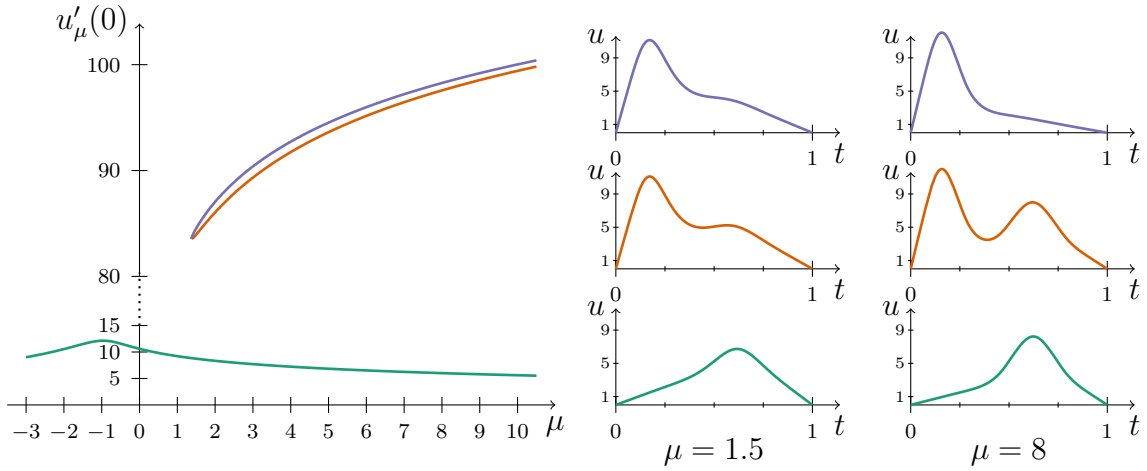
\begin{figure}[htb]
  \centering
  \begin{tikzpicture}[x=5mm, y=0.7mm]
    \draw[->] (-3.5, 0) -- (11, 0) node[below]{$\mu$};
    \draw (0, 0) -- (0, 15);
    \draw[dotted, thick] (0, 15) -- (0, 25);
    \foreach \x in {-3, -2,..., 10} {
      \draw (\x, 3pt) -- (\x, -3pt) node[below]{$\scriptstyle \x$};
    }
    \foreach \y in {5, 10, 15} {
      \draw (3pt, \y) -- (-3pt, \y) node[left]{$\scriptstyle \y$};
    }
    \ifverif
      \draw[color=c0, line width=3pt, opacity=0.3]
      plot file{/tmp/branches-4pi-0.dat};
      \draw[c0, line width=0.7pt,  mark=*, only marks, mark size=0.5pt]
      plot file{\imagepath{branches-4pi-0.dat}};
    \else
      \draw[c0, line width=1pt]
      plot file{\imagepath{branches-4pi-0.dat}};
    \fi
    \begin{scope}[yshift=-95mm, y=1.4mm]
      \draw[->] (0, 80) -- (0, 100) -- +(0, 3ex) node[left]{$u'_\mu(0)$};
      \foreach \y in {80, 90, 100} {
        \draw (3pt, \y) -- (-3pt, \y) node[left]{$\scriptstyle \y$};
      }
      \foreach \i in {1, 2}{
        \ifverif
          \draw[color=c\i, line width=3pt, opacity=0.3]
          plot file{\imagepath{branches-4pi-\i.dat}};
          \draw[c\i,  mark=*, only marks, mark size=0.5pt]
          plot file{\imagepath{branches-4pi-\i.dat}};
        \else
          \draw[c\i, line width=1pt]
          plot file{\imagepath{branches-4pi-\i.dat}};
        \fi
      }
    \end{scope}
    \foreach \x/\m in {63mm/1.5, 102mm/8}{
      \begin{scope}[xshift=\x, x=26mm, y=1.1mm]
        \node[below] at (0.5, -3mm) {$\mu = \m$};
        \graphMu[
          \foreach \x in {0.25, 0.5, 0.75}{
            \draw (\x, 1pt) -- (\x, -1pt);}
        ]{sol-4-\m}{1, 5, 9}{9.8}
      \end{scope}
    }
  \end{tikzpicture}
  \caption{The branches of positive solutions to \eqref{eq:Dmu}
    with $q(t) = \sin(4\pi t)$ in $\intervalcc{0,1}$
    and $g(u)=u^{3}$ (on the left) and
    graphs of solutions for some $\mu$ (on the right).}
  \label{fig:bifurc-4pi}
\end{figure}

\begin{figure}
  \centering
  \begin{tikzpicture}[x=6mm, y=0.75mm]
    \draw[->] (-1.5, 0) -- (10.8, 0) node[below]{$\mu$};
    \draw[->] (0, 0) -- (0, 70) -- +(0, 3ex) node[left]{$u'_\mu(0)$};
    \foreach \x in {-1,..., 10} {
      \draw (\x, 3pt) -- (\x, -3pt) node[below]{$\scriptstyle \x$};
    }
    \foreach \y in {10, 20,..., 70} {
      \draw (3pt, \y) -- (-3pt, \y) node[left]{$\scriptstyle \y$};
    }
    \begin{scope}[line width=1pt]
      \clip (-1.5, 0) rectangle (11, 73);
      \foreach \i in {0, 1, 2}{
        \draw[color=c\i] plot file{\imagepath{branches-e0.01-\i.dat}};
      }
      \node[below] at (4, 51) {$\epsilon = 0.01$};
      \foreach \i in {0, 1, 2}{
        \draw[dashed, color=c\i] plot file{\imagepath{branches-e0.1-\i.dat}};
      }
      \node[below] at (7, 67) {$\epsilon = 0.1$};
    \end{scope}
    \foreach \e/\xs in {0.01/73mm, 0.1/103mm}{
      \begin{scope}[xshift=\xs, x=21mm, y=1.1mm]
        \node[below] at (0.5, -3mm) {$\mu = 10$};
        \node[above] at (0.5, 54mm) {$\epsilon = \e$};
        \foreach \i/\ys in {0/0mm, 1/21mm, 2/42mm}{
          \begin{scope}[yshift=\ys]
            \draw[->] (0,0) -- (1, 0) -- +(2ex, 0) node[below]{$t$};
            \draw[->] (0,0) -- (0, 10) -- +(0, 1.1ex) node[left]{$u$};
            \foreach \x in {0, 1} {
              \draw (\x, 3pt) -- (\x, -3pt) node[below]{$\scriptstyle \x$};
            }
            \foreach \y in {1, 5, 9} {
              \draw (1pt, \y) -- (-1pt, \y)
              node[left]{$\scriptscriptstyle \y$};
            }
            \draw[color=c\i, line width=1pt]
            plot file{\imagepath{sol-e\e-10-\i.dat}};
          \end{scope}
        }
      \end{scope}
    }
  \end{tikzpicture}
  \caption{The branches of positive solutions to
    \eqref{eq:Dmu} with
    $h(t) = \sin(3\pi t/(1-\epsilon))$ in
    $\intervalcc{0,1-\epsilon}$, $h(t) = 0$ in
    $\intervalcc{1-\epsilon, 1}$ and $g(u)=u^{3}$ (on the left)
    and the graphs of
    solutions for $\mu = 10$ (on the right).}
  \label{fig:bifurc-epsilon}
\end{figure}
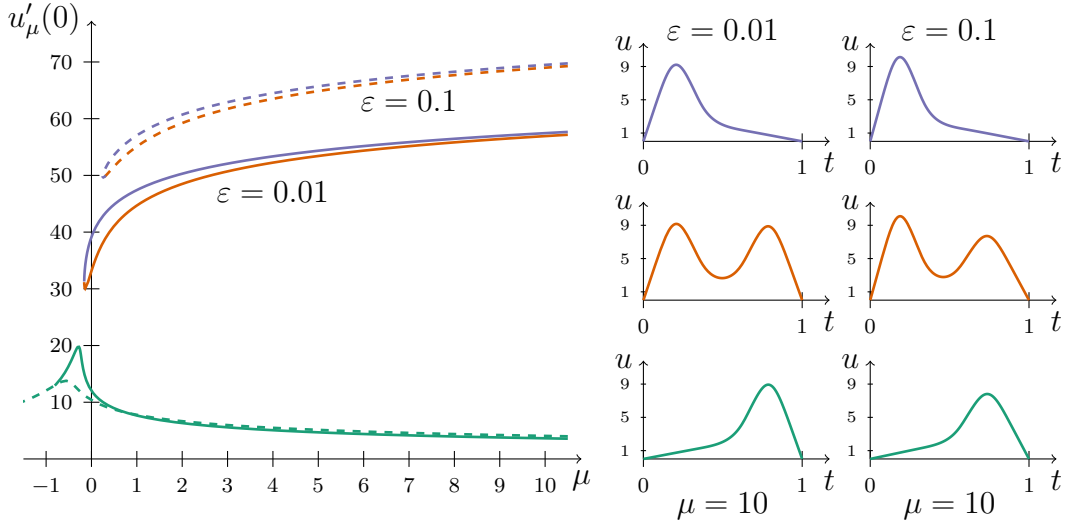

Notice that
this single bifurcation point is a rare occurrence and is not a
consequence solely of the symmetry of the weight $h$.  Indeed, on
Fig.~\ref{fig:bifurc-5pi}, six of the seven branches of positive
solutions for the symmetric weight $h(t) = \sin(5\pi t)$ collapse two
by two.  Moreover, all degenerate turning points occur at positive values
of $\mu$ (those occurring at the same $\mu$ are due to the fact that
the solutions on these branches are symmetric to each other).
The fact that all the branches but one collapse two by two also occurs
in asymmetric cases such as  $h(t) = \sin(4\pi t)$ (see
Fig.~\ref{fig:bifurc-4pi}) and for small perturbations of
$h(t) = \sin(3\pi t)$ (see Fig.~\ref{fig:bifurc-epsilon}).

\begin{figure}[htb]
  \centering
  \begin{tikzpicture}[x=26mm, y=1.2mm]
    \foreach \x/\m in {0mm/-5, 37mm/-1, 74mm/0, 111mm/1}{
      \begin{scope}[xshift=\x]
        \node[below] at (0.5, -3mm) {$\mu = \m$};
        \graphMu[
          \foreach \x in {0.25, 0.5, 0.75}{
            \draw (\x, 1pt) -- (\x, -1pt);}
        ]{sol-4-\m}{1, 5}{6}
      \end{scope}
     }
  \end{tikzpicture}
  \caption{Graphs of the unique positive solution to \eqref{eq:Dmu}
    with $q(t) = \sin(4\pi t)$ in $\intervalcc{0,1}$
    and $g(u)=u^{3}$ along the lower branch of solutions.}
  \label{fig:bifurc-4pi-neg-branch}
\end{figure}
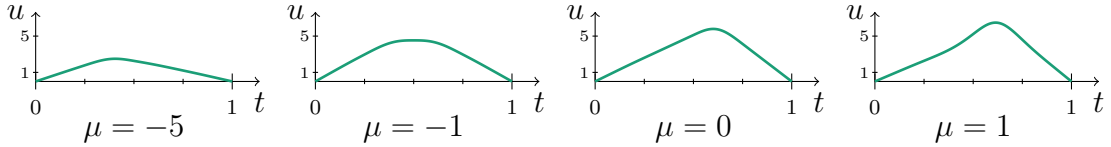

\medskip

In all the above examples (Figs.~\ref{fig:bifurc-3pi},
\ref{fig:bifurc-5pi}, \ref{fig:bifurc-4pi},
and~\ref{fig:bifurc-epsilon}), one branch of positive solutions exists
for all values of $\mu$.  This is expected as the
existence of at least one positive solution for all $\mu \in \IR$ has
been established~\cite[Theorem~5.1]{FeZa-15jde}
(extending several earlier results
\cite{de-Figueiredo-1982, ErWa-94, Nussbaum-1975} solely dealing with the
case $\mu < 0$).
Moreover,
numerical experiments show that the ground states of $\mathcal{A}$ form
such a branch and that this branch is made of unimodal solutions
(i.e., increasing then decreasing functions as they are non-negative).
When $h$ is symmetric (i.e., is even with respect to the center of the interval
$\intervaloo{a,b}$), the symmetric solutions with the lower action
$\mathcal{A}$ also form a branch living for all $\mu$ but that branch
may not coincide with the previous one as a symmetry breaking may
occur (see Fig.~\ref{fig:bifurc-3pi}), in which case  functions
along that branch may not be unimodal for all $\mu \ge 0$.

\medskip

The fact that other branches extend up to $\mu = 0$ is a delicate
question.  Clearly this is linked to the question of the uniqueness of
the positive solutions for $\mu = 0$.  In this case, the weight in
\eqref{eq:Dmu} is positive in each $I^+_i$ and identically zero in
every interval $I_{i}^{-}$.
The multiplicity of positive solutions of \eqref{eq:Dmu}
depends on the length of the intervals $I_{i}^{-}$: if the intervals
$I_{i}^{-}$ are sufficiently small, uniqueness holds (see, for
instance,
\cite{CLGT-pp, Kajikiya-Sim-Tanaka2018,
  MoNe-59}
).

\medskip

When $\mu$ is negative enough, Figs.~\ref{fig:bifurc-3pi},
\ref{fig:bifurc-5pi}, \ref{fig:bifurc-4pi},
and~\ref{fig:bifurc-epsilon} indicate that the problem admits a single
solution.
For $\mu < 0$, $h^+ - \mu h^- \ge 0$ but this uniqueness is
not a consequence of the known criteria (neither the one presented in
Section~\ref{section-3}, nor those found in the literature, see for
example~\cite{Co-67, Kw-90, Kw-91, NiNu-85, Ta-10}).  This solution is
trivially unimodal as the solutions are concave for $\mu \le 0$ (and
is symmetric when $h$ is).
As far as we know, the question of the uniqueness of positive solutions
for large negative $\mu$ is open.

\medskip

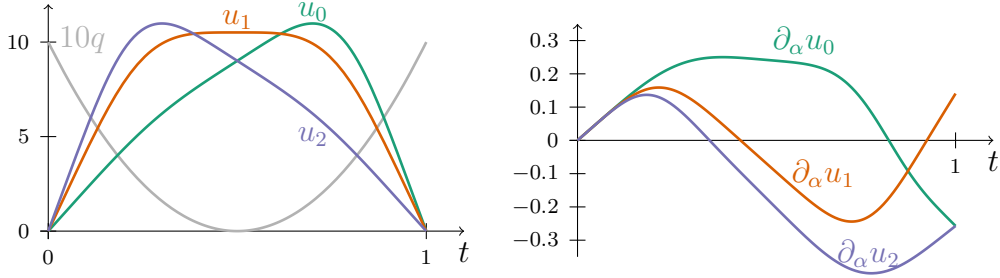
\begin{figure}[hbt]
  \centering
  \begin{tikzpicture}[x=50mm, y=2.5mm]
    \draw[->] (-3pt, 0) -- (1.1, 0) node[below]{$t$};
    \draw[->] (0, -3pt) -- (0, 12);
    \foreach \x in {0, 1} {
      \draw (\x, 3pt) -- (\x, -3pt) node[below]{$\scriptstyle \x$};
    }
    \foreach \y in {0, 5, 10} {
      \draw (3pt, \y) -- (-3pt, \y) node[left]{$\scriptstyle \y$};
    }
    \begin{scope}[color=black!30]
      \draw[line width=1pt] plot file {\imagepath{non-uniqueness-q.dat}};
      \node[right] at (0., 10.1) {$10q$};
    \end{scope}
    \foreach \i/\p in {0/{0.7, 11.6}, 1/{0.5, 11.2}, 2/{0.7, 5}} {
      \draw[color=c\i, line width=1pt]
      plot file{\imagepath{non-uniqueness-\i.dat}};
      \node[color=c\i] at (\p) {$u_{\i}$};
    }
    \begin{scope}[xshift=70mm, yshift=12mm, y=44mm]
      \draw[->] (-3pt, 0) -- (1.1, 0) node[below]{$t$};
      \draw[->] (0, -0.35) -- (0, 0.35);
      \foreach \x in {1} {
        \draw (\x, 3pt) -- (\x, -3pt) node[below]{$\scriptstyle \x$};
      }
      \foreach \y in {-0.3, -0.2, -0.1, 0, 0.1, 0.2, 0.3} {
        \draw (3pt, \y) -- (-3pt, \y) node[left]{$\scriptstyle \y$};
      }
      \foreach \i/\p in {0/{0.6, 0.3}, 1/{0.65, -0.1}, 2/{0.77, -0.34}} {
        \draw[color=c\i, line width=1pt] plot file{%
          \imagepath{non-uniqueness-\i-lin.dat}};
        \node[color=c\i] at (\p) {$\partial_\alpha u_{\i}$};
      }
    \end{scope}
  \end{tikzpicture}
  \caption{The graphs of the three positive solutions $u_i$ to the
    Dirichlet problem \eqref{eq:edo-main} when $q(t) = (2t-1)^{2}$,
    $t\in \intervalcc{0,1}$, and $g(u)=u^{3}$,
    $u\in\intervalco{0,+\infty}$ (on the left).  The
    linearization with respect to the initial velocity (on the right).}
  \label{fig:non-uniqueness}
\end{figure}

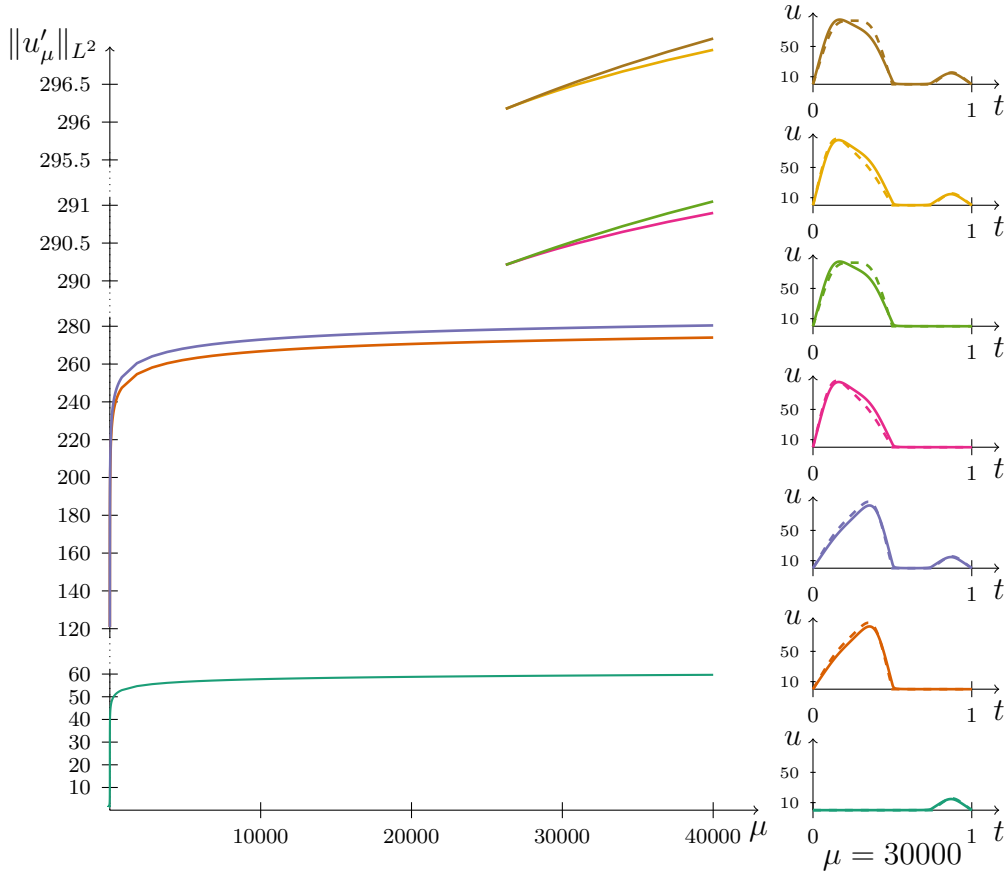
\begin{figure}
  \centering
  \begin{tikzpicture}[x=0.002mm, y=0.3mm]
    \draw[->] (-3.5, 0) -- (43000, 0) node[below]{$\mu$};
    \foreach \x in {10000, 20000, 30000, 40000} {
      \draw (\x, 3pt) -- (\x, -3pt) node[below]{$\scriptstyle \x$};
    }
    \draw[dotted] (0,0) -- (0, 290); 
    \draw (0, 0) -- (0, 62);
    \foreach \y in {10, 20,..., 60} {
      \draw (3pt, \y) -- (-3pt, \y) node[left]{$\scriptstyle \y$};
    }
    \draw[color=c0, thick] plot file{\imagepath{branches-non-uniqueness-0.dat}};
    \begin{scope}[y=0.25mm, yshift=-6mm]
      \draw (0, 115) -- (0, 285);
      \foreach \y in {120, 140,..., 280} {
        \draw (3pt, \y) -- (-3pt, \y) node[left]{$\scriptstyle \y$};
      }
      \foreach \i in {1, 2}{
        \draw[color=c\i, line width=1pt] plot
        file{\imagepath{branches-non-uniqueness-\i.dat}};
      }  
    \end{scope}
    \begin{scope}[y=10mm, yshift=-2830mm]
      \draw (0, 289.9) -- (0, 291.1);
      \foreach \y in {290, 290.5, 291} {
        \draw (3pt, \y) -- (-3pt, \y) node[left]{$\scriptstyle \y$};
      }
      \foreach \i in {3, 4}{
        \draw[color=c\i, line width=1pt] plot
        file{\imagepath{branches-non-uniqueness-\i.dat}};
      }
    \end{scope}
    \begin{scope}[y=10mm, yshift=-2869mm]
      \draw[->] (0, 295.4) -- (0, 297) node[left]{$\norm{u'_\mu}_{L^2}$};
      \foreach \y in {295.5, 296, 296.5} {
        \draw (3pt, \y) -- (-3pt, \y) node[left]{$\scriptstyle \y$};
      }
      \foreach \i in {5, 6}{
        \draw[color=c\i, line width=1pt] plot
        file{\imagepath{branches-non-uniqueness-\i.dat}};
      }
    \end{scope}
    \begin{scope}[xshift=93mm, x=21mm, y=0.1mm]
      \foreach \i/\ys in {0/0mm, 1/16mm, 2/32mm, 3/48mm, 4/64mm,
        5/80mm, 6/96mm}{
        \begin{scope}[yshift=\ys]
          \draw[->] (0,0) -- (1, 0) -- +(2ex, 0) node[below]{$t$};
          \draw[->] (0,0) -- (0, 75) -- +(0, 1.1ex) node[left]{$u$};
          \foreach \x in {0, 1} {
            \draw (\x, 3pt) -- (\x, -3pt) node[below]{$\scriptstyle \x$};
          }
          \foreach \y in {10, 50,..., 70} {
            \draw (1pt, \y) -- (-1pt, \y) node[left]{$\scriptscriptstyle \y$};
          }
          \draw[color=c\i, line width=1pt, dashed]
          plot file{\imagepath{sol-non-uniqueness-limit-\i.dat}};
          \draw[color=c\i, line width=1pt]
          plot file{\imagepath{sol-non-uniqueness-30000-\i.dat}};
        \end{scope}
      }
      \node[below] at (0.5, -3mm) {$\mu = 30000$};
    \end{scope}
  \end{tikzpicture}
  \caption{The branches $(\mu, u)$ of positive solutions
    when the limit problem possesses three positive solutions on the
    first interval $I_1^+$ when $h$ is given by \eqref{eq:h-3sols}
    and $g(u) = u^3$ (on the left).  The
    graphs of solutions for $\mu = 30000$ as well as those of the
    limit profiles, dashed (on the right).}
  \label{fig:branches-non-uniqueness}
\end{figure}

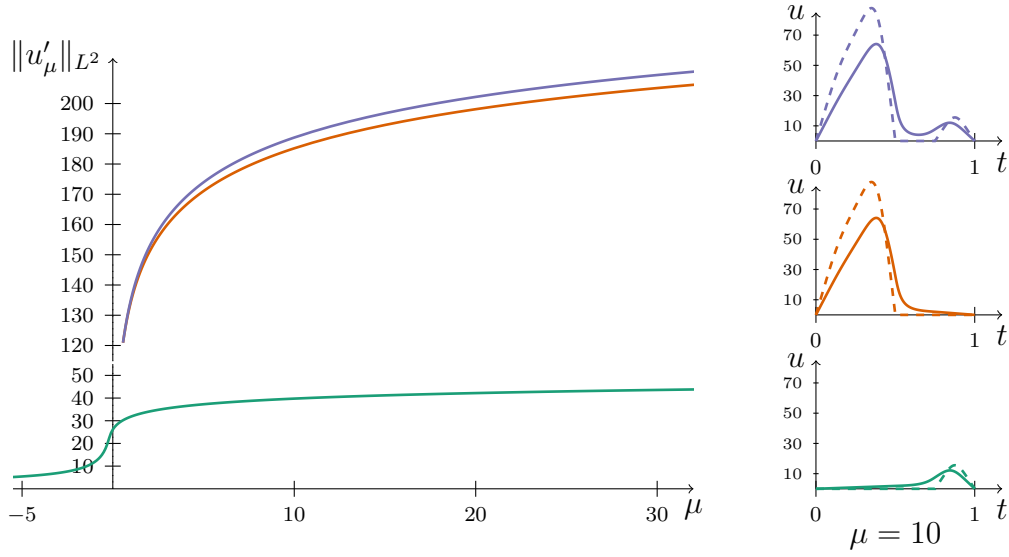
\begin{figure}[htb]
  \centering
  \begin{tikzpicture}[x=2.4mm, y=0.3mm]
    \draw[->] (-5.5, 0) -- (32, 0) node[below]{$\mu$};
    \foreach \x in {-5, 10, 20, 30} {
      \draw (\x, 3pt) -- (\x, -3pt) node[below]{$\scriptstyle \x$};
    }
    \draw[dotted] (0,0) -- (0, 100);
    \draw (0,0) -- (0, 55);
    \foreach \y in {10, 20,..., 50} {
      \draw (3pt, \y) -- (-3pt, \y) node[left]{$\scriptstyle \y$};
    }
    \begin{scope}
      \clip (-5.5, 0) rectangle (32, 50);
      \draw[color=c0, line width=1pt] plot
      file{\imagepath{branches-non-uniqueness-zoom-0.dat}};
    \end{scope}
    \begin{scope}[y= 0.4mm, yshift=-16ex]
      \draw[->] (0, 115) -- (0, 215) node[left]{$\norm{u'_\mu}_{L^2}$};
      \foreach \y in {120, 130,..., 200} {
        \draw (3pt, \y) -- (-3pt, \y) node[left]{$\scriptstyle \y$};
      }
      \clip (0, 115) rectangle (32, 215);
      \foreach \i in {1, 2}{
        \draw[color=c\i, line width=1pt] plot
        file{\imagepath{branches-non-uniqueness-zoom-\i.dat}};
      }
    \end{scope}
    \begin{scope}[xshift=93mm, x=21mm, y=0.2mm]
      \foreach \i/\ys in {0/0mm, 1/23mm, 2/46mm}{
        \begin{scope}[yshift=\ys]
          \draw[->] (0,0) -- (1, 0) -- +(2ex, 0) node[below]{$t$};
          \draw[->] (0,0) -- (0, 75) -- +(0, 1.1ex) node[left]{$u$};
          \foreach \x in {0, 1} {
            \draw (\x, 3pt) -- (\x, -3pt) node[below]{$\scriptstyle \x$};
          }
          \foreach \y in {10, 30,..., 70} {
            \draw (1pt, \y) -- (-1pt, \y) node[left]{$\scriptscriptstyle \y$};
          }
          \draw[color=c\i, line width=1pt, dashed]
          plot file{\imagepath{sol-non-uniqueness-limit-\i.dat}};
          \draw[color=c\i, line width=1pt]
          plot file{\imagepath{sol-non-uniqueness-10-\i.dat}};
        \end{scope}
      }
      \node[below] at (0.5, -3mm) {$\mu = 10$};
    \end{scope}
  \end{tikzpicture}
  \caption{The branches $(\mu, u)$ of positive solutions
    when the limit problem possesses three positive solutions on the
    first interval $I_1^+$ when $h$ is given by \eqref{eq:h-3sols}
    and $g(u) = u^3$ (on the left) and the
    graphs of solutions for $\mu = 10$ as well as those of the
    limit profiles, dashed (on the right).}
  \label{fig:branches-non-uniqueness-zoom}
\end{figure}

It is well known that there are non-negative weights $q$ such that
the Dirichlet problem \eqref{eq:edo-main} possesses several positive
solutions.  In 1959,
Moore and Nehari \cite{MoNe-59} gave such an example as a smooth
perturbation of a piece-wise constant function. 
Here we have chosen the simple $\mathcal{C}^\infty$ symmetric function
\begin{equation*}
  q\colon \intervalcc{0,1} \to \IR,\  q(t) = (2t-1)^2.
\end{equation*}
Still considering $g(u) = u^{3}$, Fig.~\ref{fig:non-uniqueness} shows
on the left the graphs of the three positive solutions $u_i$,
$i \in \{0, 1, 2\}$, to problem \eqref{eq:edo-main}.  If we denote
$u(\cdot; \alpha)$ the solution to the differential equation of
\eqref{eq:edo-main} for the initial conditions (see
Section~\ref{section-3})
\begin{equation*}
  u(0; \alpha) = 0, \qquad
  u'(0; \alpha) = \alpha,
\end{equation*}
then $u_i(\cdot) = u(\cdot; \alpha_i)$ for some $\alpha_i > 0$.  The
right graphs of Fig.~\ref{fig:non-uniqueness} show the derivatives
$t \mapsto \partial_\alpha u(t; \alpha_i)$ of the solutions with
respect to
the initial velocity $\alpha$.  Based on the numerical evidence of
these graphs, we
posit that $\partial_\alpha u(1; \alpha_i) \ne 0$ and so that
the three solutions $u_i$ are non-degenerate.  Therefore
a branch should emanate from each of them if $h$ has the same
shape as the above $q$
on an interval $I^+_i$.
More precisely, let us consider the function
\begin{equation}
  \label{eq:h-3sols}
  h \colon \intervalcc{0,1} \to \IR,\  h(t) =
  \begin{cases}
    (t - 1/4)^2& \text{if } t \in \intervalco{0, 1/2},\\
    -\sin(4\pi t)& \text{if } t \in \intervalcc{1/2, 1}
  \end{cases}
\end{equation}
for problem~\eqref{eq:Dmu}.  Since for this $h$ there are three
positive limit profiles to choose from on $I^+_1 = \intervalcc{0,1/2}$
and one positive limit
profile on $I^+_2 = \intervalcc{3/4, 1}$, we expect
$4 \cdot 2 - 1 = 7$ branches of positive solutions.  This is confirmed
by our numerical experiments, see
Figs.~\ref{fig:branches-non-uniqueness}
and~\ref{fig:branches-non-uniqueness-zoom}.  Four of theses branches
live until $\mu \approx 26250$ while two persist until
$\mu \approx 0.56$ and the last one lives for all $\mu$ (and,
according to the graphs, is made of unimodal functions).
We believe that the techniques developed in this paper can be extended
to prove the existence of these seven branches and more generally to
establish the exact multiplicity result under a sole condition of the
non-degeneracy of the positive limit profiles.

Let us conclude with a more technical remark.   Notice
that the vertical axis on Figs.~\ref{fig:branches-non-uniqueness}
and~\ref{fig:branches-non-uniqueness-zoom} reports $\norm{u'_\mu}_{L^2}$ and
not $u'_\mu(0)$.  The reason is that the initial velocities of some
solutions are really close, especially for the large values of $\mu$
we need to tackle to see all branches.  This is particularly true for
solutions having the same limit profile on $I^+_1$ but differing on
$I^+_2$.  Thus, instead of tracking $u'_\mu(0)$, we use a finite
element discretization of \eqref{eq:Dmu} and start a branch
continuation algorithm at
$\mu = 10^6$ from the limit profiles (the dashed graphs on
Figs.~\ref{fig:branches-non-uniqueness}
and~\ref{fig:branches-non-uniqueness-zoom}) refined thanks to a damped
Newton's method.

\section*{Data sharing}

Data sharing is not applicable to this article.

\section*{Acknowledgments}

During the time the research was conducted, the first author has been supported by the project ERC Advanced Grant 2013 n.~339958 ``Complex Patterns for Strongly Interacting Dynamical Systems - COMPAT'',  by the Belgian F.R.S.-FNRS - Fonds de la Recherche Scientifique, \textit{Charg\'{e} de recherches} project: ``Quantitative and qualitative properties of positive solutions to indefinite problems arising from population genetic models: topological methods and numerical analysis'', and also partially by the
Grup\-po Na\-zio\-na\-le per l'Anali\-si Ma\-te\-ma\-ti\-ca, la Pro\-ba\-bi\-li\-t\`{a} e le lo\-ro
Appli\-ca\-zio\-ni (GNAMPA) of the Isti\-tu\-to Na\-zio\-na\-le di Al\-ta Ma\-te\-ma\-ti\-ca (INdAM), in particular by the INdAM-GNAMPA project ``Analisi qualitativa di problemi differenziali non lineari''.

\IfFileExists{images.tex}{%

}{%
  \bibliographystyle{elsart-num-sort}
  \bibliography{FeTr_biblio}
}

\end{document}